\documentclass[12pt]{amsart}
\usepackage[margin=1.2in]{geometry}
\usepackage{amssymb}
\usepackage{mathtools}
\usepackage{enumitem}
\usepackage{tikz}
\usetikzlibrary{arrows.meta,calc}
\usepackage[colorlinks,linkcolor=blue,citecolor=blue,urlcolor=blue]{hyperref}

\numberwithin{equation}{section}

\theoremstyle{plain}
\newtheorem{theorem}{Theorem}[section]
\newtheorem{lemma}[theorem]{Lemma}
\newtheorem{proposition}[theorem]{Proposition}
\newtheorem{corollary}[theorem]{Corollary}
\newtheorem*{theoremA}{Theorem A}
\newtheorem*{theoremB}{Theorem B}
\newtheorem*{theoremC}{Theorem C}
\newtheorem*{theoremD}{Theorem D}
\newtheorem*{theoremE}{Theorem E}
\theoremstyle{definition}
\newtheorem{definition}[theorem]{Definition}
\newtheorem{example}[theorem]{Example}
\newtheorem{remark}[theorem]{Remark}

\newcommand{\R}{\mathbb{R}}
\newcommand{\scal}[2]{\langle #1,#2\rangle}
\newcommand{\Mreg}{M_{\mathrm{reg}}}
\newcommand{\hT}{\hat\tau}
\DeclareMathOperator{\conv}{conv}
\DeclareMathOperator{\tr}{tr}
\DeclareMathOperator{\ext}{ext}

\newcommand{\hullpath}{(-1.2,0.9) arc[start angle=90,end angle=281.06,radius=0.9]
     -- (2.015,-0.289) arc[start angle=-78.94,end angle=90,radius=0.6] -- cycle}
\newcommand{\drawsection}{\draw[sect] (-1.2,0) circle (0.9);
  \draw[sect] (1.9,0.3) circle (0.6);
  \draw[sect] (0.9,-0.28) ellipse [x radius=0.18, y radius=0.09];
  \draw[chord] (-1.2,0.9) -- (1.9,0.9);
  \draw[chord] (-1.027,-0.883) -- (2.015,-0.289);}

\begin{document}

\title{Support functions of sections of $p$-minimal tubes}

\author{Vladimir G. Tkachev}
\address{Department of Mathematics, Link\"oping University, SE-581 83 Link\"oping, Sweden}
\email{vladimir.tkatjev@liu.se}



\begin{abstract}
Let $M$ be a compact $p$-minimal hypersurface in $\R^{n+1}$ with boundary in two parallel hyperplanes. We study the convex hulls $\Omega(\tau)$ of its horizontal sections through their support function $H(\tau,\cdot)$, without any assumption on the topology of $M$. We prove that $H$ is a viscosity subsolution of a degenerate elliptic inequality. For $n=2$ it is of Monge--Amp\`ere type, and its equality case for minimal surfaces with circular sections is Riemann's system. A posteriori, $H$ is locally semiconvex and every hull is a smooth convex body; the inequality then also holds in the sense of measures, with the edges of the hulls appearing as atoms of Alexandrov's area functions.

Averaging over directions gives a catenary-type inequality for the mean width of the hulls. It yields lifetime bounds for $p$-minimal tubes when $n>p$. For minimal surfaces in $\R^3$ it gives a sharp bound for the distance between the planes in terms of the perimeters of the hulls of the boundary curves, with equality only for a catenoid. For one convex curve in each plane this recovers, by a different method, the nonexistence and rigidity theorem of Bernstein and Breiner. The bound extends to non-convex curves and to boundaries with several components. Finally, we show that the Meeks annulus conjecture is equivalent to the absence of atoms in the area functions of the hulls of all intermediate sections.
\end{abstract}

\maketitle

\section{Introduction}\label{sec:intro}

\subsection{Background}
How far apart can two parallel planes be if curves in them bound a connected minimal surface? This question has a long history. Nitsche gave a bound for minimal annuli in terms of the diameters of the boundary curves (see \cite{Rossman08}), and Rossman \cite{Rossman08} extended a bound of this type to surfaces of arbitrary topology. For annuli, Osserman and Schiffer \cite{OssermanSchiffer75} found the convexity of the lengths of the sections in conformal coordinates that governs the problem. Building on it, Bernstein and Breiner \cite{BernsteinBreiner14} obtained the sharp answer in terms of the lengths of the boundary curves: the threshold is set by marginally stable catenoids, it holds for immersed surfaces of any genus, and equality characterizes the catenoid. At the level of sections, Shiffman \cite{Shiffman56} showed that every intermediate section of a minimal annulus bounded by convex curves is convex, and a circle if the boundary curves are circles. Riemann \cite{Riemann1867} and Enneper \cite{Enneper1869} had classified the minimal surfaces foliated by circles in parallel planes.

In higher dimensions the corresponding objects are minimal \emph{tubes}: hypersurfaces $M^n\subset\R^{n+1}$ whose boundary lies in two parallel hyperplanes. For $n\ge3$, Miklyukov and Vedenyapin \cite{MV86} showed that their \emph{lifetime}, the distance between the hyperplanes, is bounded in terms of the size of the sections; see also \cite{MT89}, and Klyachin \cite{Klyachin92} for arbitrary codimension. In \cite{T97} the author extended this to $p$-minimal tubes, by comparing $M$ with the circumscribed balls of its sections.

A further motivation is the \emph{Meeks annulus conjecture} (MAC) \cite{Meeks78}: a compact connected minimal surface in $\R^3$ bounded by two convex Jordan curves in parallel planes is an annulus. Such surfaces are embedded by Ekholm, White and Wienholtz \cite{EWW02}. MAC was proved by Schoen \cite{Schoen83} for boundaries with two planes of symmetry, by Meeks and White \cite{MeeksWhite91,MeeksWhite93} for stable surfaces, by Ros \cite{Ros96} when one curve is a vertical translate of the other, using the L\'opez--Ros deformation, and by Traizet \cite{Traizet10} for planes close enough to each other, for surfaces of genus one, and for any fixed genus when the two curves are symmetric about a common line \cite[Cor.~1, Cor.~2]{Traizet10}. It is open in general, even for two circles with different axes \cite[\S1]{Traizet10}.

In this paper we show that the convex hulls of the sections of such surfaces obey a comparison principle that needs no information on their topology. From it we obtain, in a unified way:
\begin{itemize}
\item lifetime bounds for $p$-minimal tubes;
\item a sharp bound for the distance between the planes in terms of the perimeters of the hulls of the boundary curves, attained only by a catenoid. This recovers the theorem of Bernstein and Breiner by a different method, and extends it to non-convex curves and to boundaries with several components;
\item a reformulation of MAC in terms of the convex geometry of the sections.
\end{itemize}

\subsection{Hulls of sections}
Let $\Sigma(\tau)\subset\R^n$ be the section of $M$ at height $\tau$, i.e. the projection of $M\cap\{x_{n+1}=\tau\}$ to $\R^n$. In \cite{T97} the circumscribed balls of the sections of a $p$-minimal hypersurface $M^n\subset\R^{n+1}$ (see Section~\ref{subsec:mainresults} for the definition) were compared with $M$ at contact points, which gives
\begin{equation}\label{eq:T97}
RR''\ge\beta(1+R'^2)+\min\{\beta,1\}\,|\xi'|^2,\qquad \beta=\tfrac{n-1}{p-1},
\end{equation}
for the radius $R(\tau)$ and centre $\xi(\tau)$ of the ball \cite[Thm.~4]{T97}. In an unpublished preprint \cite{T92} the balls were replaced by the convex hulls
\begin{equation}\label{eq:SigmaOmega}
\Omega(\tau)=\conv\Sigma(\tau),\qquad\text{where}\quad \Sigma(\tau)=\{y\in\R^n\mid (y,\tau)\in M\},
\end{equation}
with second derivatives controlled through difference-of-convex functions. The text of \cite{T92} is not available to us, and the exposition below is self-contained. We work with the support function \cite[\S1.7]{Schneider14}, \cite{Leichtweiss80}
\begin{equation}\label{eq:Hdef}
H(\tau,w)=\max_{y\in\Omega(\tau)}\scal{y}{w},\qquad w\in\R^n .
\end{equation}

The hulls see much more than the balls: the shape of the sections, their edges and their perimeters. Two observations make them usable.
\begin{itemize}
\item \emph{The boundary of a hull sees only regular points of the height} (Lemma~\ref{lem:regular}). So the comparison never meets the critical points of $x_{n+1}$, and no information on the topology of $M$ is needed.
\item \emph{Regularity comes for free.} The comparison is first proved in the viscosity sense \cite{CIL92}, which requires no regularity of $H$.
It then forces $H$ to be semiconvex, i.e.\ convex up to a quadratic term, which recovers the difference-of-convex structure of \cite{T92} as a consequence rather than an assumption.
\end{itemize}

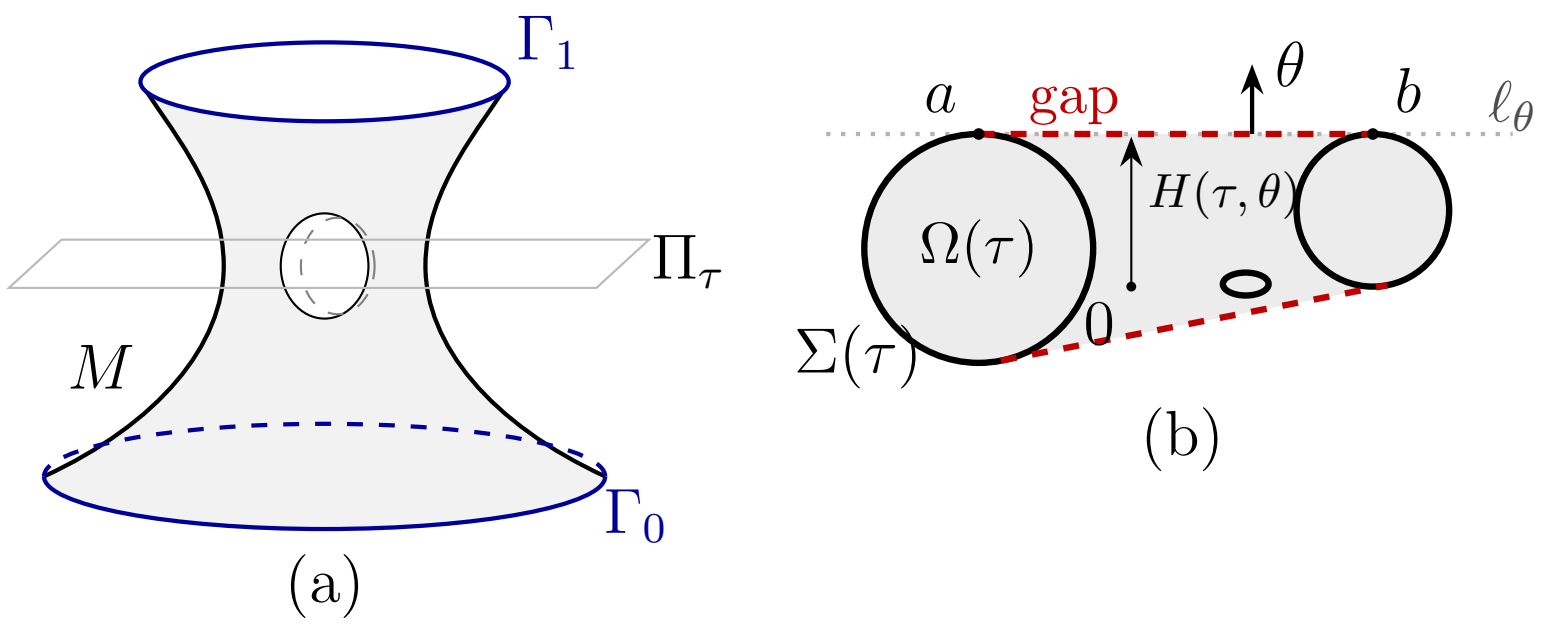
\begin{figure}[ht]
\centering
\begin{tikzpicture}[scale=0.58,
  sect/.style={very thick, black},
  hull/.style={fill=gray!15},
  chord/.style={very thick, dashed, red!75!black},
  bnd/.style={thick, blue!60!black}]
\begin{scope}[yshift=-0.2cm]
 \fill[gray!10] (-3.2,0) .. controls (-1.9,0.6) and (-1.15,1.5) .. (-1.15,2.4)
     .. controls (-1.15,3.3) and (-1.8,4.0) .. (-2.1,4.5)
     -- (2.1,4.5) .. controls (1.8,4.0) and (1.15,3.3) .. (1.15,2.4)
     .. controls (1.15,1.5) and (1.9,0.6) .. (3.2,0)
     arc[start angle=0,end angle=-180,x radius=3.2,y radius=0.6] -- cycle;
 \draw[thick] (-3.2,0) .. controls (-1.9,0.6) and (-1.15,1.5) .. (-1.15,2.4)
     .. controls (-1.15,3.3) and (-1.8,4.0) .. (-2.1,4.5);
 \draw[thick] (3.2,0) .. controls (1.9,0.6) and (1.15,1.5) .. (1.15,2.4)
     .. controls (1.15,3.3) and (1.8,4.0) .. (2.1,4.5);
 \draw[bnd] (3.2,0) arc[start angle=0,end angle=-180,x radius=3.2,y radius=0.6];
 \draw[bnd, dashed] (3.2,0) arc[start angle=0,end angle=180,x radius=3.2,y radius=0.6];
 \fill[white] (0,4.5) ellipse [x radius=2.1, y radius=0.45];
 \draw[bnd] (0,4.5) ellipse [x radius=2.1, y radius=0.45];
 \fill[white] (0,2.4) ellipse [x radius=0.5, y radius=0.6];
 \draw[thin] (0,2.4) ellipse [x radius=0.5, y radius=0.6];
 \draw[thin, dashed, gray] (0.15,2.4) ellipse [x radius=0.42, y radius=0.55];
 \draw[gray!55] (-3.6,2.15) -- (3.1,2.15) -- (3.7,2.7) -- (-3.0,2.7) -- cycle;
 \node[blue!60!black] at (3.55,-0.45) {$\Gamma_0$};
 \node[blue!60!black] at (2.55,4.95) {$\Gamma_1$};
 \node at (-2.55,1.25) {$M$};
 \node at (4.15,2.45) {$\Pi_\tau$};
 \node at (0,-1.25) {(a)};
\end{scope}
\begin{scope}[xshift=9.2cm, yshift=2.4cm, scale=1.45]
 \fill[hull] \hullpath;
 \draw[gray!60, dotted, thick] (-2.4,0.9) -- (3.0,0.9);
 \node[gray!60!black, font=\small] at (3.0,1.12) {$\ell_\theta$};
 \drawsection
 \fill (-1.2,0.9) circle (1.3pt) node[above left] {$a$};
 \fill (1.9,0.9) circle (1.3pt) node[above right] {$b$};
 \fill (0.0,-0.3) circle (1.1pt) node[below left=-1pt] {$0$};
 \draw[-{Stealth[length=2mm]}] (0.0,-0.3) -- (0.0,0.88);
 \node[font=\scriptsize, right=-1pt] at (0.0,0.42) {$H(\tau,\theta)$};
 \draw[-{Stealth[length=2mm]}, thick] (0.95,0.9) -- (0.95,1.45) node[right] {$\theta$};
 \node at (-2.15,-0.85) {$\Sigma(\tau)$};
 \node[font=\small] at (-1.2,0) {$\Omega(\tau)$};
 \node[red!75!black, font=\small] at (-0.45,1.12) {gap};
 \node at (0.4,-1.5) {(b)};
\end{scope}
\end{tikzpicture}
\caption{(a) A schematic, hypothetical connected surface $M$ of genus one spanning convex curves $\Gamma_0,\Gamma_1$. The front and back walls are joined around a window (thin contour; its hidden mirror image is dashed), and the plane $\Pi_\tau$ through the window cuts $M$ in a disconnected section. (b) The plane $\Pi_\tau$: the section $\Sigma(\tau)$ (solid), its hull $\Omega(\tau)$ (shaded) with support function $H(\tau,\theta)$, and two \emph{supporting chords} (Definition~\ref{def:bridge}): edges of $\Omega(\tau)$ with endpoints on $\Sigma(\tau)$ whose interiors (dashed) miss $\Sigma(\tau)$. The endpoints $a,b$ are regular points of $x_3$ (Lemma~\ref{lem:regular}) and, for a.e.\ $\tau$, points of tangency of $M$ with one plane (Proposition~\ref{prop:bridge}). The small inner component is invisible to $H$. By Theorem~E, for convex $\Gamma_0,\Gamma_1$, $M$ is an annulus if and only if no section has a chord.}
\label{fig:bridge}
\end{figure}

\subsection{Main results}\label{subsec:mainresults}
Let $p>1$, $\Pi_\tau=\{x_{n+1}=\tau\}$ and $\tau_0<\tau_1$. A $C^2$-immersed hypersurface $M^n\subset\R^{n+1}$ is \emph{$p$-minimal} if its mean curvature $H_M$ and the normal curvature $k_e$ in the direction of steepest ascent of the height satisfy $H_M=-(p-2)k_e$ wherever $\nabla x_{n+1}\ne0$ (Section~\ref{sec:prelim}). Then the height $x_{n+1}$ is $p$-harmonic in the induced metric \cite[Def.~2]{T97}; for $p=2$, $M$ is minimal. Write $x=(y,\tau)\in\R^n\times\R$.

\begin{theoremA}[Theorem~\ref{thm:main}]
Let $M^n\subset\R^{n+1}$ be a compact connected $C^2$-immersed $p$-minimal hypersurface with $\partial M\subset\Pi_{\tau_0}\cup\Pi_{\tau_1}$. For $\tau\in[\tau_0,\tau_1]$ let $\Sigma(\tau)$ be its section at height $\tau$, $\Omega(\tau)=\conv\Sigma(\tau)$, and $H(\tau,\cdot)$ the support function of $\Omega(\tau)$. Then $H$ is a viscosity subsolution, in the sense of Definition~\ref{def:visc}, of
\begin{equation}\label{eq:main}
H_{\tau\tau}\ \ge\ \frac{1+H_\tau^2}{p-1}\,\tr\mathcal R^{-1}+\scal{\mathcal R^{-1}b}{b},
\end{equation}
where $\mathcal R=D^2_wH(\tau,\theta)|_{\theta^\perp}$ is the matrix of principal radii of curvature of $\partial\Omega(\tau)$ in the outer normal direction $\theta$, and $b$ is the projection of $\nabla_wH_\tau(\tau,\theta)$ onto $\theta^\perp$.
\end{theoremA}

When $\Omega(\tau)$ is a ball, \eqref{eq:main} reduces to the inequality behind \eqref{eq:T97} (Corollary~\ref{cor:balls}). The proof rests on the fact that every point of $\partial\Omega(\tau)\cap\Sigma(\tau)$ is a regular point of the height function (Lemma~\ref{lem:regular}). The same first-order argument shows that the hulls have no corners.

\begin{theoremB}[Lemma~\ref{lem:smooth}, Theorem~\ref{thm:semiconvex}, Corollary~\ref{cor:pointwise}]
Let $M$ and $H$ be as in Theorem~A. Then for every $\tau\in(\tau_0,\tau_1)$, $\Omega(\tau)$ is a smooth convex body. The function $H$ is locally semiconvex on $(\tau_0,\tau_1)\times(\R^n\setminus0)$, hence locally a difference of convex functions. In particular, $H$ is twice differentiable almost everywhere, \eqref{eq:main} holds pointwise almost everywhere, and the singular part of the Hessian of $H$ is positive semidefinite.
\end{theoremB}

For surfaces in $\R^3$ we write
\begin{equation}\label{eq:hdef}
h(\tau,\varphi)=H\big(\tau,(\cos\varphi,\sin\varphi)\big).
\end{equation}

\begin{theoremC}[Theorem~\ref{thm:MA}, Propositions~\ref{prop:equality} and~\ref{prop:MAmeasure}]
Let $M$ and $H$ be as in Theorem~A with $n=2$, and let $h$ be given by \eqref{eq:hdef}. Then $h$ is a viscosity subsolution of the Monge--Amp\`ere type inequality
\begin{equation}\label{eq:MA}
h_{\tau\tau}\,(h_{\varphi\varphi}+h)-h_{\tau\varphi}^2\ \ge\ \frac{1+h_\tau^2}{p-1},\qquad h_{\varphi\varphi}+h>0 .
\end{equation}
The matrix on the left is a nonnegative matrix-valued measure whose absolutely continuous part satisfies \eqref{eq:MA}, and whose second diagonal entry is the family of area functions $S(\Omega(\tau),\cdot)$. Equality holds for surfaces foliated by simple closed curves of positive curvature. For $p=2$ and circular sections of radius $R(\tau)$ and centre $\xi(\tau)$, equality is Riemann's system
\begin{equation}\label{eq:RiemannIntro}
RR''=1+R'^2+|\xi'|^2,\qquad R\xi''=2R'\xi' .
\end{equation}
\end{theoremC}

Averaging \eqref{eq:main} over the sphere of directions gives a closed differential inequality for the \emph{mean half-width} of the hulls,
\begin{equation}\label{eq:Bdef}
B(\tau)=\frac{1}{|S^{n-1}|}\int_{S^{n-1}}H(\tau,\theta)\,d\theta .
\end{equation}
With $\beta=\frac{n-1}{p-1}$, it reads
\begin{equation}\label{eq:Cateq}
B''\ \ge\ \beta\,\frac{1+B'^2}{B-\lambda},
\end{equation}
where $\lambda(\tau)\ge0$ measures the flat part of the boundary of the hull $\Omega(\tau)$: it is the normalized mass of the singular part of the measure $\Delta_Sh+(n-1)h$ on $S^{n-1}$, a multiple of the first area measure of $\Omega(\tau)$ (see \ref{S6} and Proposition~\ref{prop:Bineq}). It vanishes if $\partial\Omega(\tau)$ is $C^2$ with positive curvature, and for $n=2$, $2\pi\lambda(\tau)$ is at least the total length of the edges of $\Omega(\tau)$.
The equality case for $n=p=2$ in \eqref{eq:Cateq} is the catenoid.

\begin{theoremD}[Theorems~\ref{thm:integrated} and~\ref{thm:rigidity}, Corollaries~\ref{cor:BB} and~\ref{cor:strictseveral}]
The following hold for compact connected hypersurfaces with boundary in $\Pi_{\tau_0}\cup\Pi_{\tau_1}$.
\begin{enumerate}[label=(\roman*)]
\item Let $M$ be as in Theorem~A with $n>p$, let $\beta=\frac{n-1}{p-1}$, and let $B$ be given by \eqref{eq:Bdef}. Then
\[
\tau_1-\tau_0\le2\kappa_\beta\,B_{\min},\qquad \kappa_\beta=\int_1^\infty\frac{dx}{\sqrt{x^{2\beta}-1}},\quad B_{\min}=\min_\tau B(\tau).
\]
\item Let $M\subset\R^3$ be a compact connected immersed minimal surface with $\partial M=\Gamma_0\cup\Gamma_1$, where each $\Gamma_i\subset\Pi_{\tau_i}$ is a finite union of pairwise disjoint closed curves, and let $\rho_i=\operatorname{Per}(\conv\Gamma_i)/(2\pi)$. Then $\tau_1-\tau_0\le D(\rho_0,\rho_1)$, where
\[
D(\rho_0,\rho_1)=\max_{0<m\le\min\{\rho_0,\rho_1\}}m\Big[\operatorname{arccosh}\frac{\rho_0}{m}+\operatorname{arccosh}\frac{\rho_1}{m}\Big]
\]
is the largest height of a piece of a catenoid bounded by coaxial circles of radii $\rho_0,\rho_1$. The inequality is strict if the hulls of the sections have edges on a set of levels of positive measure, in particular if $M$ is smooth up to the boundary and no component of $\Gamma_0$, or no component of $\Gamma_1$, has the same convex hull as the whole of $\Gamma_0$, respectively $\Gamma_1$.
\item Let $M$ be as in (ii), and assume moreover that $x$ is an immersion up to the boundary and that $x|_{\partial M}$ is injective. If $\tau_1-\tau_0=D(\rho_0,\rho_1)$, then $\Gamma_0,\Gamma_1$ are coaxial circles and $M$ is the extremal catenoid piece between them.
\end{enumerate}
\end{theoremD}

Part (i) sharpens \cite{MV86,MT89} and \cite[Cor.~5]{T97}: the circumradius is replaced by the smaller quantity $B_{\min}$, and no condition on the sections is needed. Parts (ii) and (iii) are compared with the literature in Section~\ref{subsec:previous}.

Finally, we apply the method to MAC. A \emph{supporting chord} of a section, briefly a \emph{chord}, is a segment $[a,b]\subset\partial\Omega(\tau)$ with $a,b\in\Sigma(\tau)$ and $(a,b)\cap\Sigma(\tau)=\emptyset$ (Definition~\ref{def:bridge}); see Figure~\ref{fig:bridge}.

\begin{theoremE}[Theorem~\ref{thm:MAC}, Corollary~\ref{cor:C1}]
Let $M\subset\R^3$ be a compact connected embedded minimal surface with $\partial M=\Gamma_0\cup\Gamma_1$, where $\Gamma_i\subset\Pi_{\tau_i}$ are Jordan curves. If
\begin{equation}\label{eq:F}
\partial\Omega(\tau)\subset\Sigma(\tau)\qquad\text{for every }\tau\in(\tau_0,\tau_1),
\end{equation}
then $M$ is an annulus and every intermediate section is a convex Jordan curve. If $\Gamma_0,\Gamma_1$ are convex, then conversely every annulus satisfies \eqref{eq:F}. Moreover, the following are equivalent:
\begin{itemize}
\item \eqref{eq:F} holds;
\item $H(\tau,\cdot)\in C^1(\R^2\setminus0)$ for every $\tau\in(\tau_0,\tau_1)$;
\item no intermediate section has a chord;
\item the area functions $S(\Omega(\tau),\cdot)$ have no atoms.
\end{itemize}
\end{theoremE}

Thus MAC reduces to the exclusion of chords. By Theorems~B and~E, the hulls of the sections of a hypothetical counterexample are smooth convex bodies that fail to be strictly convex. The first part of Theorem~E requires no convexity of $\Gamma_i$. The equivalences rest on the fact that a section of a minimal surface contains no straight segment (Lemma~\ref{lem:nosegment}). Theorem~E thus recasts MAC as a statement about the convex geometry of the sections: a counterexample must produce atoms in the area functions of the hulls of its sections. In this form the obstruction is local and measure-theoretic, and comparison arguments of the type developed here apply to it.

\subsection{Relation to previous work}\label{subsec:previous}
\emph{Height bounds in $\R^3$.} For one convex curve in each plane, Theorem~D(ii)--(iii) recovers the nonexistence and rigidity theorem of Bernstein and Breiner \cite[Thm.~4.1, Prop.~4.4]{BernsteinBreiner14}: the two thresholds coincide (Section~\ref{sec:integrated}). The methods are different.
\begin{itemize}
\item Bernstein and Breiner use the conformal structure of annuli and the Osserman--Schiffer inequality, and extend the argument to arbitrary genus. They also prove area-minimizing properties of the catenoid, which we do not address.
\item Our argument uses only a maximum principle at the boundary of the hulls. It is insensitive to the conformal type, and it controls the hull of \emph{every} intermediate section, not only the boundary curves. This is what makes the extensions possible: lengths are replaced by hull perimeters, which is a strict improvement for non-convex curves; boundaries with several components are allowed; and the bound is strict as soon as the hulls have edges.
\end{itemize}
In terms of diameters, the bound improves Rossman's inequality $\tau_1-\tau_0\le\max_i\operatorname{diam}\Gamma_i$ \cite{Rossman08} to $\tau_1-\tau_0\le\max_i\operatorname{diam}\Gamma_i/\sinh x_0\approx0.6627\,\max_i\operatorname{diam}\Gamma_i$. Here $x_0\approx1.1997$ is the positive root of $x\tanh x=1$, which determines the marginally stable catenoid (Section~\ref{sec:integrated}).

\emph{Lifetime of tubes.} For $n>p$, Theorem~D(i) contains the lifetime bounds of \cite{MV86,MT89} for $p=2$ and of \cite[Cor.~5]{T97}, with a smaller constant and without the condition on the sections used in \cite{T97}. We do not treat higher codimension, where the methods of Klyachin \cite{Klyachin92} would be needed (Remark~\ref{rem:regular}).

\emph{Sections of minimal surfaces.} Shiffman's theorem \cite{Shiffman56} starts from an annulus and uses a global conformal parametrization in which the height is a coordinate. The hull method uses no topology, and Theorem~C recovers the equation of Riemann's examples as its equality case. In this sense, Theorem~A may be viewed as a convex-hull version of Shiffman's convexity, valid for surfaces of arbitrary topology, where convexity of the sections themselves is not available.

In a companion paper \cite{T26notes} we study the mechanism by which chords are born, using the Weierstrass representation, and the number of components of the sections.

\subsection*{Organization}
Section~\ref{sec:prelim} collects preliminaries on $p$-minimal hypersurfaces and support functions. Section~\ref{sec:contact} proves that contact points are regular and that the hulls are smooth. Section~\ref{sec:main} proves Theorem~A and the semiconvexity of $H$, and Section~\ref{sec:MA} treats the case $n=2$. Section~\ref{sec:integrated} contains the integrated inequality, its consequences and the comparison with known results. Section~\ref{sec:MAC} applies the method to MAC, and Section~\ref{sec:several} treats boundaries with several components.

\section{Preliminaries}\label{sec:prelim}

Throughout Sections~\ref{sec:prelim}--\ref{sec:integrated}, $p>1$ and $M^n\subset\R^{n+1}$ is a compact connected $C^2$-immersed $p$-minimal hypersurface, $x\colon M\to\R^{n+1}$ is the immersion, and $x(\partial M)\subset\Pi_{\tau_0}\cup\Pi_{\tau_1}$ with $\tau_0<\tau_1$. The sections $\Sigma(\tau)$, their hulls $\Omega(\tau)$ and the support function $H$ are defined by \eqref{eq:SigmaOmega} and \eqref{eq:Hdef}.

Let $e=e_{n+1}$. For $m\in M$ let $\nu(m)$ be a unit normal and $A$ the Weingarten map, so that $\nabla^2 x_j(X,Y)=\scal{AX}{Y}\scal{\nu}{e_j}$ and $\Delta x=H_M\nu$, where $H_M=\tr A$ \cite[Ch.~VII]{Osserman86}. The gradient of a linear function $\scal{x}{v}$ on $M$ is the tangential projection $v^\top$; in particular $\nabla x_{n+1}=e^\top$. Let
\[
\Mreg=\{m\in M^\circ:\ e^\top(m)\ne0\},\qquad \hT=e^\top/|e^\top|\ \text{ on }\Mreg ,
\]
and $k_e=\scal{A\hT}{\hT}$, the normal curvature in the direction of steepest ascent of $x_{n+1}$. As in \eqref{eq:SigmaOmega}, the sections and their convex hulls are
\[
\Sigma(\tau)=\{y\in\R^n\mid (y,\tau)\in x(M)\},\qquad \Omega(\tau)=\conv\Sigma(\tau).
\]

Since $M$ is connected and $x_{n+1}$ takes the values $\tau_0$ and $\tau_1$ on $\partial M$, every $\Sigma(\tau)$, $\tau\in[\tau_0,\tau_1]$, is nonempty and compact. By the maximum principle for the convex hull \cite[Cor.~3]{T97}, $x(M)\subset\conv x(\partial M)$, so $x_{n+1}(M)\subset[\tau_0,\tau_1]$.

We say that $M$ is \emph{$p$-minimal} if
\begin{equation}\label{eq:pminimal}
H_M=-(p-2)k_e\qquad\text{on }\Mreg .
\end{equation}
By the computation in the proof of \cite[Prop.~1]{T97}, on $\Mreg$
\[
\Delta_px_{n+1}=\operatorname{div}\big(|\nabla x_{n+1}|^{p-2}\nabla x_{n+1}\big)=|\nabla x_{n+1}|^{p-2}\scal{\nu}{e}\big(H_M+(p-2)k_e\big).
\]
Hence the height of a $p$-minimal hypersurface is $p$-harmonic. Conversely, if $x_{n+1}$ is $p$-harmonic, then \eqref{eq:pminimal} holds wherever $\scal{\nu}{e}\ne0$, hence on the closure of this set in $\Mreg$. The proviso cannot be dropped: the height on a vertical cylinder $\Gamma\times\R$ over a closed curve $\Gamma\subset\R^n$ is $p$-harmonic for every $p$, including $p=2$, while $k_e=0$ and $H_M\ne0$ there. Thus $p$-harmonicity of the height is equivalent to \eqref{eq:pminimal} if the set of points of $\Mreg$ with vertical tangent hyperplane has empty interior. The converse direction of \cite[Prop.~1]{T97} must be read with this assumption.
On $\Mreg$ define the symmetric bilinear form and the operator
\begin{equation}\label{eq:G}
G(X,Y)=\scal{X}{Y}+(p-2)\scal{X}{\hT}\scal{Y}{\hT},\qquad Lu=\Delta u+(p-2)\nabla^2u(\hT,\hT).
\end{equation}
Since $\hT=\nabla x_{n+1}/|\nabla x_{n+1}|$, we have $\Delta_px_{n+1}=|\nabla x_{n+1}|^{p-2}Lx_{n+1}$ on $\Mreg$. This splitting of the $p$-Laplacian into the Laplacian and the second derivative in the gradient direction is also the starting point of Moser's approximation of the inverse mean curvature flow by $p$-harmonic functions \cite{Moser07}, where it is read on the level sets of a $p$-harmonic function; see \cite{MariRigoliSetti22,BenattiMariRigoliSettiXu25} for recent developments. Here it is read on the hypersurface $M$ itself, and \eqref{eq:pminimal} is the resulting curvature identity.
The eigenvalues of $G$ are $1$ (multiplicity $n-1$) and $p-1$, so $L$ is uniformly elliptic on $\Mreg$, and for $F\in C^2(\R^k)$ and $u_1,\dots,u_k\in C^2$,
\begin{equation}\label{eq:chain}
L\,F(u_1,\dots,u_k)=\sum_i F_i\,Lu_i+\sum_{i,j}F_{ij}\,G(\nabla u_i,\nabla u_j).
\end{equation}

\begin{lemma}\label{lem:L}
On $\Mreg$ every coordinate function is $L$-harmonic: $Lx_j=0$, $1\le j\le n+1$.
\end{lemma}

\begin{proof}
By the formulas above,
\[Lx_j=\big(H_M+(p-2)\scal{A\hT}{\hT}\big)\scal{\nu}{e_j}=\big(H_M+(p-2)k_e\big)\scal{\nu}{e_j},
\]
which vanishes by \eqref{eq:pminimal}.
\end{proof}

Thus \eqref{eq:pminimal} is equivalent to the $L$-harmonicity of \emph{all} coordinates, for an operator depending only on the direction $\hT$. For $p=2$, $L=\Delta$, and Lemma~\ref{lem:L} is the classical characterization of minimal immersions \cite[\S4]{Osserman86}. Harmonicity of the height alone does not imply minimality, as the vertical cylinders above show.

\section{Support functions}
We collect standard facts \cite[\S\S1.4, 1.7, 2.5]{Schneider14}, \cite{Leichtweiss80}. Let $K\subset\R^n$ be compact and $\Omega=\conv K$, with support function $H_\Omega(w)=\max_{y\in\Omega}\scal{y}{w}$.
\begin{enumerate}[label=(S\arabic*)]
\item\label{S1} $H_\Omega$ is convex and positively $1$-homogeneous, and $H_\Omega=H_K$.
\item\label{S2} $H_{t\Omega_1\oplus s\Omega_2}=tH_{\Omega_1}+sH_{\Omega_2}$ for $t,s\ge0$; $\Omega_1\subset\Omega_2$ iff $H_{\Omega_1}\le H_{\Omega_2}$.
\item\label{S3} $H_\Omega$ is differentiable at $\theta\ne0$ iff the face $F(\Omega,\theta)=\{y\in\Omega:\scal{y}{\theta}=H_\Omega(\theta)\}$ is a single point $y$, and then $\nabla H_\Omega(\theta)=y$.
\item\label{S4} If $\partial\Omega$ is of class $C^2$ with positive curvature, then $H_\Omega\in C^2(\R^n\setminus 0)$ and the eigenvalues of $D^2H_\Omega(\theta)|_{\theta^\perp}$ are the principal radii of curvature of $\partial\Omega$ at $\nabla H_\Omega(\theta)$. For $n=2$, writing $h(\varphi)=H_\Omega(\cos\varphi,\sin\varphi)$, the radius of curvature is $h''+h$.
\item\label{S5} $\ext\Omega\subset K$, i.e. every extreme point of $\Omega$ belongs to $K$.
\item\label{S6} \emph{Area and curvature measures.} Let $h=H_\Omega|_{S^{n-1}}$. The distribution
\[
\sigma(\Omega,\cdot)=\Delta_Sh+(n-1)h
\]
on $S^{n-1}$, where $\Delta_S$ is the Laplacian of $S^{n-1}$, is a nonnegative Radon measure of total mass $(n-1)\int_{S^{n-1}}h$. It equals $(n-1)S_1(\Omega,\cdot)$, where $S_1$ is the first area measure \cite[\S4.2]{Schneider14}. We write 
\[
\sigma(\Omega,\cdot)=\sigma^{ac}(\Omega,\cdot)+\sigma^{s}(\Omega,\cdot)
\] 
for its \emph{Lebesgue decomposition} with respect to the standard measure $d\theta$ on $S^{n-1}$; more generally, for a Radon measure $\mu$ we denote by $\mu^{ac}$ and $\mu^{s}$ its absolutely continuous and singular parts with respect to the relevant Lebesgue or spherical measure. If $\partial\Omega$ is $C^2$ with positive curvature, the density of $\sigma(\Omega,\cdot)$ is $\tr D^2H_\Omega|_{\theta^\perp}$, the sum of the principal radii.
\end{enumerate}

For $n=2$, $\sigma(\Omega,\cdot)=h''+h$ is Alexandrov's \emph{area function}, or length measure, $S(\Omega,\cdot)$. For a Borel set $\eta\subset S^1$, $S(\Omega,\eta)$ is the length of the set of boundary points having an outer normal in $\eta$. Its total mass is the perimeter of $\Omega$, and its atoms are the edges: $S(\Omega,\{\theta\})$ is the length of the face $F(\Omega,\theta)$.

Dually, Alexandrov's \emph{integral curvature} $\omega(\Omega,\beta)$ of a Borel set $\beta\subset\partial\Omega$ is the measure of the set of outer unit normals of $\Omega$ at points of $\beta$ \cite{Alexandrov06}, \cite[\S4.2]{Schneider14}. For $n=2$ its atoms are the corners, and $\omega(\Omega,\{y\})$ is the angle of the normal cone at $y$.

The following is \cite[Thm.~2]{T97} rewritten via \ref{S2}.

\begin{lemma}\label{lem:convex}
For every $w\in\R^n$ the function $\tau\mapsto H(\tau,w)$ is convex on $[\tau_0,\tau_1]$. In particular it is locally Lipschitz on $(\tau_0,\tau_1)$, and $H$ is continuous on $(\tau_0,\tau_1)\times\R^n$.
\end{lemma}

\begin{proof}
By \cite[Thm.~2]{T97}, $\Omega(t\tau_1'+(1-t)\tau_2')\subset t\,\Omega(\tau_1')\oplus(1-t)\,\Omega(\tau_2')$, which by \ref{S2} is the convexity of $H(\cdot,w)$. Joint continuity follows from
\[|H(\tau,w)-H(\tau,w')|\le \max_{M}|y|\,|w-w'|.
\]
\end{proof}

\section{Contact points are regular}\label{sec:contact}

\begin{definition}\label{def:contact}
Let $\tau\in(\tau_0,\tau_1)$ and $|\theta|=1$. A point $m\in M$ is a \emph{contact point} for $(\tau,\theta)$ if $x(m)=(y,\tau)$ with $\scal{y}{\theta}=H(\tau,\theta)$.
\end{definition}

Contact points exist since $\Sigma(\tau)$ is compact, and every point of $\partial\Omega(\tau)\cap\Sigma(\tau)$ is (the image of) a contact point for some $\theta$. Since $\partial M$ lies in $\Pi_{\tau_0}\cup\Pi_{\tau_1}$, contact points are interior points of $M$.

The basic inequality of the method is
\begin{equation}\label{eq:basic}
\scal{y(m')}{w}\le H(\tau(m'),w)\qquad\text{for all } m'\in M,\ w\in\R^n,
\end{equation}
where $x(m')=(y(m'),\tau(m'))$. It follows from $y(m')\in\Sigma(\tau(m'))\subset\Omega(\tau(m'))$.

\begin{lemma}\label{lem:regular}
Every contact point $m$ lies in $\Mreg$, i.e. $\nabla x_{n+1}(m)\ne0$.
\end{lemma}

\begin{proof}
Suppose $e^\top(m)=0$. Then $T_mM=e^\perp$, so $(\theta,0)\in T_mM$. By \eqref{eq:basic}, the function
\[
\psi(m')=\scal{y(m')}{\theta}-H(\tau(m'),\theta)
\]
satisfies $\psi\le0=\psi(m)$ near $m$. Since $\nabla\tau(m)=0$ and $x\in C^2$, we have $|\tau(m')-\tau|=O(d(m',m)^2)$. Since $H(\cdot,\theta)$ is Lipschitz near $\tau$ (Lemma~\ref{lem:convex}), the second term of $\psi$ also changes by $O(d(m',m)^2)$. The first term has gradient $(\theta,0)^\top=(\theta,0)\ne0$ at $m$. Hence $\psi$ takes positive values along a curve through $m$ in the direction $(\theta,0)$, a contradiction.
\end{proof}

\begin{remark}\label{rem:regular}
Lemma~\ref{lem:regular} uses neither convexity of $\partial M$ nor any information on the topology of $M$. It does use the codimension one: for higher codimension $T_mM$ may be orthogonal to both $e$ and $(\theta,0)$.
At this point, the methods of Klyachin \cite{Klyachin92} would be needed.
\end{remark}

The same first-order argument shows that the hulls have no corners.

\begin{lemma}\label{lem:smooth}
Let $\tau\in(\tau_0,\tau_1)$. Then $\Omega(\tau)$ has interior points, and at every point of $\partial\Omega(\tau)$ it has a unique outer unit normal. In other words, $\Omega(\tau)$ is a smooth convex body \cite[\S2.2]{Schneider14}. More precisely, let $\theta$ be an outer unit normal of $\Omega(\tau)$ at $a\in\partial\Omega(\tau)\cap\Sigma(\tau)$, and let $m$ be a point with $x(m)=(a,\tau)$. Then
\begin{equation}\label{eq:normalc}
\nu(m)\parallel(\theta,-c)\qquad\text{for some }c\in\big[\partial_\tau^-H(\tau,\theta),\,\partial_\tau^+H(\tau,\theta)\big].
\end{equation}
\end{lemma}

\begin{proof}
The proof adapts a standard convex-geometric argument: one first identifies the normal at a regular contact point, then excludes lower-dimensional affine degeneracy, and finally uses the face structure of a convex set to prove uniqueness of the supporting normal.

\emph{Step 1: the normal at a contact point.} The point $m$ is a contact point for $(\tau,\theta)$, hence regular by Lemma~\ref{lem:regular}. Let $g=H(\cdot,\theta)$, which is convex and locally Lipschitz (Lemma~\ref{lem:convex}). By \eqref{eq:basic},
\[
\psi(m')=\scal{y(m')}{\theta}-g(\tau(m'))\ \le\ 0=\psi(m).
\]
Let $X\in T_mM$ and let $\gamma$ be a $C^2$ curve with $\gamma(0)=m$, $\gamma'(0)=X$.
\begin{itemize}
\item If $d\tau(X)=0$, then $\tau(\gamma(s))-\tau=O(s^2)$, hence $g(\tau(\gamma(s)))-g(\tau)=O(s^2)$, and $\scal{y(\gamma(s))-a}{\theta}\le O(s^2)$ for both signs of $s$. Thus $d\scal{y}{\theta}(X)=0$.
\item Consequently $(\theta,0)^\top$ is orthogonal to $\ker d\tau$, i.e. $(\theta,0)^\top=c\,e^\top$ for some $c\in\R$.
\item If $d\tau(X)>0$, then for small $s>0$
\[
c\,d\tau(X)\,s+o(s)=\scal{y(\gamma(s))-a}{\theta}\le g(\tau(\gamma(s)))-g(\tau)=\partial^+_\tau g\;d\tau(X)\,s+o(s),
\]
so $c\le\partial^+_\tau g$. Replacing $X$ by $-X$ gives $c\ge\partial_\tau^-g$.
\end{itemize}
Hence $(\theta,-c)^\top=(\theta,0)^\top-c\,e^\top=0$, which is \eqref{eq:normalc}.

\emph{Step 2: interior points.} Suppose $\Omega(\tau)$ lies in an affine hyperplane $A$ of $\R^n$, and let $a$ be an extreme point of $\Omega(\tau)$. By \ref{S5}, $a\in\Sigma(\tau)$; since every point of $\Omega(\tau)$ is a contact point for the normal of $A$, any $m$ over $a$ is regular. Hence $x$ maps a neighbourhood of $m$ in $x_{n+1}^{-1}(\tau)$ onto an embedded $(n-1)$-dimensional $C^1$ submanifold through $(a,\tau)$. It lies in $\Sigma(\tau)\times\{\tau\}\subset A\times\{\tau\}$, so by invariance of domain it contains a neighbourhood of $a$ in $A$. Thus $a$ is the midpoint of a segment in $\Omega(\tau)$, which contradicts extremality.

\emph{Step 3: unique normals.} Let $a\in\partial\Omega(\tau)\cap\Sigma(\tau)$, and let $N(a)$ be the normal cone. Fix one $m$ over $a$. By Step 1, every unit $\theta\in N(a)$ is $\pm\nu_h(m)/|\nu_h(m)|$, where $\nu_h\ne0$ is the horizontal part of $\nu$. By Step 2, $N(a)$ contains no line. Hence $N(a)$ is a single ray.

Now let $z\in\partial\Omega(\tau)\setminus\Sigma(\tau)$. The point $z$ lies in the relative interior of a unique face $F$ of $\Omega(\tau)$ \cite[\S2.1]{Schneider14}. Let $a$ be an extreme point of $F$; it is extreme in $\Omega(\tau)$, so $a\in\Sigma(\tau)$ by \ref{S5}. If $\theta\in N(z)$, the linear function $\scal{\cdot}{\theta}$ attains its maximum over $F$ at the relative interior point $z$, so it is constant on $F$, and $\theta\in N(a)$. Thus $N(z)\subset N(a)$ is a single ray.
\end{proof}

\begin{remark}\label{rem:smooth}
Lemma~\ref{lem:smooth} uses only Lemma~\ref{lem:regular}, so it holds for immersed hypersurfaces and every $p>1$. In the language of \ref{S6} for $n=2$: the integral curvature $\omega(\Omega(\tau),\cdot)$ has no atoms, and the area function $S(\Omega(\tau),\cdot)$ charges every open arc of $S^1$. The dual properties, that $S(\Omega(\tau),\cdot)$ has no atoms and that $\omega(\Omega(\tau),\cdot)$ charges every open arc of $\partial\Omega(\tau)$, both mean that $\Omega(\tau)$ has no edges. For embedded minimal surfaces in $\R^3$ bounded by convex curves, the absence of edges at all levels is equivalent to MAC (Corollary~\ref{cor:C1}, Theorem~\ref{thm:MAC}).
\end{remark}

\section{The main inequality}\label{sec:main}

\begin{definition}\label{def:visc}
Let $U$ be continuous on $(\tau_0,\tau_1)\times(\R^n\setminus0)$ and $1$-homogeneous in $w$. We say that $U$ satisfies \eqref{eq:main} \emph{in the viscosity sense} if the following holds at every $(\tau,\theta)$, $|\theta|=1$. Let $P$ be any $C^2$ function near $(\tau,\theta)$, $1$-homogeneous in $w$, with $P\ge U$ near $(\tau,\theta)$ and $P(\tau,\theta)=U(\tau,\theta)$. Then $\mathcal R_P=D^2_wP(\tau,\theta)|_{\theta^\perp}$ is positive definite, and with $b_P=\mathrm{pr}_{\theta^\perp}\nabla_wP_\tau(\tau,\theta)$,
\begin{equation}\label{eq:mainP}
P_{\tau\tau}\ \ge\ \frac{1+P_\tau^2}{p-1}\,\tr\mathcal R_P^{-1}+\scal{\mathcal R_P^{-1}b_P}{b_P}.
\end{equation}
\end{definition}

This is the usual notion of viscosity subsolution \cite[Def.~2.2]{CIL92}, \cite[Def.~1.2.1]{NTV14}: the inequality is tested on $C^2$ functions touching $U$ from above. The operator in \eqref{eq:main} is degenerate elliptic in the Hessian; for viscosity theory of degenerate elliptic equations see \cite[\S1.5]{NTV14} and the references there. The homogeneity of $P$ is no restriction: if $P\ge U$ touches at $(\tau,\theta)$, so does $|w|\,P(\tau,w/|w|)$.

\begin{lemma}\label{lem:algebra}
Let $s\in\R$, $\theta\in S^{n-1}$, $N=(\theta,-s)\in\R^{n+1}$, $\nu=N/|N|$, and let ${}^\top$ denote the orthogonal projection onto $\nu^\perp$. For a unit $u\perp\theta$ and $\beta\in\R$ put $V=(u,-\beta)$. Then
\begin{equation}\label{eq:three}
|e^\top|^2=\frac{1}{1+s^2},\qquad |V^\top|^2=1+\frac{\beta^2}{1+s^2},\qquad \scal{V^\top}{e^\top}=-\frac{\beta}{1+s^2}.
\end{equation}
\end{lemma}

\begin{proof}
We have $\scal{e}{N}=-s$, $\scal{V}{N}=\beta s$ and $|N|^2=1+s^2$. Hence $|e^\top|^2=1-s^2/(1+s^2)$,
\[
|V^\top|^2=1+\beta^2-\beta^2s^2/(1+s^2),
\]
and
\[\scal{V^\top}{e^\top}=\scal{V}{e}-\scal{V}{\nu}\scal{e}{\nu}=-\beta+\beta s^2/(1+s^2)
\]
implying \eqref{eq:three}.
\end{proof}

\begin{theorem}\label{thm:main}
The support function $H$ satisfies \eqref{eq:main} in the viscosity sense on $(\tau_0,\tau_1)\times(\R^n\setminus 0)$. In particular, \eqref{eq:main} holds pointwise at every $(\tau,\theta)$ near which $H$ is of class $C^2$.
\end{theorem}

\begin{proof}
The argument proceeds by first-order variations in $w$ and on $M$, which identify the contact point and its normal, followed by a second-order variation. Diagonalizing the resulting quadratic inequality with respect to $\mathcal R_P$ and optimizing it gives both $\mathcal R_P>0$ and \eqref{eq:mainP}. Let $P$ touch $H$ from above at $(\tau,\theta)$ as in Definition~\ref{def:visc}, and let $m$ be a contact point for $(\tau,\theta)$.

\emph{Step 1: first order in $w$.} By \eqref{eq:basic}, $\Psi(m',w')=\scal{y(m')}{w'}-P(\tau(m'),w')\le0$ for $(m',w')$ near $(m,\theta)$, with equality at $(m,\theta)$. Since $w'$ ranges over an open set, $\nabla_w\Psi(m,\theta)=0$:
\begin{equation}\label{eq:y=gradP}
y(m)=\nabla_wP(\tau,\theta).
\end{equation}
In particular the face $F(\Omega(\tau),\theta)$ is the single point $y(m)$ by \ref{S3}. By Lemma~\ref{lem:regular}, $m\in\Mreg$.

\emph{Step 2: first order on $M$.} Let $w:M\to\R^n$ be any $C^2$ map defined near $m$ with $w(m)=\theta$, and put $\psi(m')=\Psi(m',w(m'))$. Then $\psi\le0=\psi(m)$, so $\nabla\psi(m)=0$ and, since $L$ is elliptic, $L\psi(m)\le0$. Using $\nabla y_i=e_i^\top$, $\nabla\tau=e^\top$ and \eqref{eq:y=gradP},
\[
\nabla\psi(m)=\sum_i\theta_ie_i^\top-P_\tau e^\top=N^\top,\qquad N:=(\theta,-P_\tau(\tau,\theta)).
\]
Hence $N\perp T_mM$, and we may take $\nu(m)=N/|N|$.

\emph{Step 3: second order.} By Lemma~\ref{lem:L}, $Ly_i=L\tau=0$. By \eqref{eq:chain} and \eqref{eq:y=gradP}, the terms containing $Lw_i$ cancel at $m$, and with $W_i=\nabla w_i(m)$,
\begin{equation}\label{eq:Phi}
0\ge L\psi(m)=2\sum_iG(a_i,W_i)-\sum_{i,j}P_{w_iw_j}\,G(W_i,W_j)-P_{\tau\tau}\,G(e^\top,e^\top),
\end{equation}
where $a_i=(e_i-P_{\tau w_i}e)^\top$ and all derivatives of $P$ are taken at $(\tau,\theta)$. Choosing $w_i=\theta_i+\ell_i$ with $\ell_i(m)=0$, the vectors $W_1,\dots,W_n\in T_mM$ are arbitrary.

\emph{Step 4: diagonalization.} Let $\{\theta,u_1,\dots,u_{n-1}\}$ be an orthonormal basis of $\R^n$ in which $\mathcal R_P$ is diagonal, $\mathcal R_Pu_\alpha=\rho_\alpha u_\alpha$. Since $P$ is $1$-homogeneous, $D^2_wP\,\theta=0$ and $\sum_i\theta_iP_{\tau w_i}=P_\tau$. Since $H(\tau,\cdot)$ is convex and $P(\tau,\cdot)$ touches it from above at $\theta$, the form $D_w^2P(\tau,\theta)$ is positive semidefinite, so $\rho_\alpha\ge0$. Rewriting \eqref{eq:Phi} in the new basis, the $\theta$-component of $a$ is $\sum_i\theta_ia_i=N^\top=0$, and with $\tilde W_\alpha=\sum_i u_{\alpha,i}W_i$ we get
\begin{equation}\label{eq:Phi2}
2\sum_{\alpha}G(\tilde a_\alpha,\tilde W_\alpha)-\sum_\alpha\rho_\alpha\,G(\tilde W_\alpha,\tilde W_\alpha)\le P_{\tau\tau}\,G(e^\top,e^\top)
\end{equation}
for all $\tilde W_\alpha\in T_mM$, where $\tilde a_\alpha=V_\alpha^\top$, $V_\alpha=(u_\alpha,-b_\alpha)$ and $b_\alpha=\scal{\nabla_wP_\tau}{u_\alpha}$.

\emph{Step 5: evaluation.} By Lemma~\ref{lem:algebra} with $s=P_\tau$, $\beta=b_\alpha$ and \eqref{eq:G},
\begin{align*}
G(e^\top,e^\top)&=(p-1)|e^\top|^2=\frac{p-1}{1+P_\tau^2},\\
G(\tilde a_\alpha,\tilde a_\alpha)&=|V_\alpha^\top|^2+(p-2)\frac{\scal{V_\alpha^\top}{e^\top}^2}{|e^\top|^2}=1+\frac{(p-1)\,b_\alpha^2}{1+P_\tau^2}.
\end{align*}
If some $\rho_\alpha=0$, take $\tilde W_\alpha=t\tilde a_\alpha$ and the other $\tilde W$'s zero; the left side of \eqref{eq:Phi2} is
\[
2t\,G(\tilde a_\alpha,\tilde a_\alpha)\ge 2t\to\infty,
\]
a contradiction. Hence $\mathcal R_P>0$. Maximizing each term of \eqref{eq:Phi2} at $\tilde W_\alpha=\tilde a_\alpha/\rho_\alpha$ gives
\[
\sum_\alpha\frac{1}{\rho_\alpha}\Big(1+\frac{(p-1)b_\alpha^2}{1+P_\tau^2}\Big)\le\frac{(p-1)P_{\tau\tau}}{1+P_\tau^2},
\]
which is \eqref{eq:mainP}, since $\sum_\alpha\rho_\alpha^{-1}=\tr\mathcal R_P^{-1}$ and $\sum_\alpha b_\alpha^2/\rho_\alpha=\scal{\mathcal R_P^{-1}b_P}{b_P}$.
\end{proof}

\begin{remark}\label{rem:meaning}
When $H$ is smooth, \eqref{eq:y=gradP} identifies the contact point as $y(\tau,\theta)=\nabla_wH(\tau,\theta)$, $\mathcal R$ is the matrix of principal radii of $\partial\Omega(\tau)$ at $y$ by \ref{S4}, and $b=\mathrm{pr}_{\theta^\perp}\partial_\tau y$ is the tangential velocity of the contact point as $\tau$ varies. The proof shows that $M$ cannot be touched from outside, at a contact point, by a hypersurface whose sections are convex with radii $\mathcal R$ and which violates \eqref{eq:main}.
\end{remark}

\subsection{The case of balls}\label{subsec:balls}
Theorem~\ref{thm:main} contains the comparison with circumscribed balls of \cite{T97}. Let $B(\tau)$ be the circumscribed ball of $\Sigma(\tau)$, with centre $\xi(\tau)$ and radius $R(\tau)$. By \ref{S2},
\begin{equation}\label{eq:PB}
P_B(\tau,w)=\scal{\xi(\tau)}{w}+R(\tau)|w|\ \ge\ H(\tau,w),
\end{equation}
with equality at $(\tau,\theta)$ whenever $\xi(\tau)+R(\tau)\theta\in\Sigma(\tau)$. The functions $R$ and $\xi$ need not be smooth, and we use them through second-order expansions. Recall that a function $f$ is \emph{twice differentiable at $\tau$ in the sense of Peano} if
\begin{equation}\label{eq:Peano}
f(\tau+t)=f(\tau)+f'(\tau)\,t+\tfrac12 f''(\tau)\,t^2+o(t^2)\qquad (t\to0).
\end{equation}
This is the sense in which convex functions, and differences of convex functions, are twice differentiable almost everywhere by Alexandrov's theorem \cite[\S6.4]{EvansGariepy15}.

\begin{corollary}\label{cor:balls}
Let $\tau\in(\tau_0,\tau_1)$, and suppose that $R$ and $\xi$ are twice differentiable at $\tau$ in the sense of Peano. Then for every unit vector $\theta$ with $\xi(\tau)+R(\tau)\theta\in\Sigma(\tau)$,
\begin{equation}\label{eq:balls}
RR''-\beta(1+R'^2)\ \ge\ (\beta-1)\scal{\xi'}{\theta}^2+|\xi'|^2+\scal{\theta}{2\beta R'\xi'-R\xi''},\qquad \beta=\tfrac{n-1}{p-1},
\end{equation}
where $R,\xi$ and their derivatives are evaluated at $\tau$.
\end{corollary}

\begin{proof}
Let $R_2,\xi_2$ be the second-order Taylor polynomials of $R,\xi$ at $\tau$, and for $\varepsilon>0$ put
\[
P_\varepsilon(\tau+t,w)=\scal{\xi_2(\tau+t)}{w}+\big(R_2(\tau+t)+\varepsilon t^2\big)|w| .
\]
This function is smooth and $1$-homogeneous in $w$ near $(\tau,\theta)$. By \eqref{eq:Peano}, the differences $\xi-\xi_2$ and $R-R_2$ are $o(t^2)$, so $P_\varepsilon\ge P_B\ge H$ near $(\tau,\theta)$ by \eqref{eq:PB}, with equality at $(\tau,\theta)$. At $(\tau,\theta)$ we have
\[
P_\tau=R'+\scal{\xi'}{\theta},\qquad P_{\tau\tau}=R''+\scal{\xi''}{\theta}+2\varepsilon,\qquad \mathcal R_P=R\,\mathrm{Id},\qquad b_P=\xi'-\scal{\xi'}{\theta}\theta .
\]
Hence \eqref{eq:mainP} gives
\[
R\big(R''+\scal{\xi''}{\theta}+2\varepsilon\big)\ \ge\ \beta\big(1+(R'+\scal{\xi'}{\theta})^2\big)+|\xi'|^2-\scal{\xi'}{\theta}^2 .
\]
Letting $\varepsilon\to0$ and rearranging gives \eqref{eq:balls}.
\end{proof}

Inequality \eqref{eq:balls} is exactly inequality \cite[(q)]{T97} in the proof of \cite[Thm.~4]{T97}, from which \eqref{eq:T97} follows by choosing $\theta$ via \cite[Lemma~3]{T97}. Thus the inequality behind \cite[Thm.~4]{T97} is a special case of Theorem~\ref{thm:main} at every $\tau$ where $R$ and $\xi$ are twice differentiable in the sense of Peano. In \cite{T97} this differentiability almost everywhere was provided by the condition $\sigma(\Sigma(\tau))\ge\varepsilon$, which makes $\xi$ a difference of convex functions. Corollary~\ref{cor:balls} does not reprove \eqref{eq:T97} without such a condition. What the present approach removes, in Section~\ref{sec:integrated}, is the need for it: there we work with $H$ instead of $(R,\xi)$.

\subsection{Semiconvexity of the support function}
The function $H$ is convex in $w$ by \ref{S1} and in $\tau$ by Lemma~\ref{lem:convex}, but separate convexity gives no control of the mixed derivatives. We now show that Theorem~\ref{thm:main} itself supplies this control: $H$ is locally semiconvex in $(\tau,w)$.

Recall that a function $f$ on an open set $U\subset\R^N$ is \emph{semiconvex} with constant $C\ge0$ if $f+\frac C2|x|^2$ is convex on every convex subset of $U$; equivalently, if $f$ is continuous and
\[
f(x+h)+f(x-h)-2f(x)\ \ge\ -C|h|^2\qquad\text{whenever }[x-h,x+h]\subset U .
\]
It is \emph{locally semiconvex} if every point of $U$ has a convex neighbourhood on which $f$ is semiconvex for some constant. A locally semiconvex function is locally Lipschitz, and locally a difference of two convex functions. By Alexandrov's theorem it is twice differentiable almost everywhere in the sense of \eqref{eq:Peano}, see \cite[\S6.4]{EvansGariepy15} and \cite[Ch.~2]{CannarsaSinestrari04}.

We use the following viscosity characterization of convexity; cf.\ \cite{Oberman07}, where convexity is characterized by the equation $\lambda_{\min}(D^2v)=0$.

\begin{lemma}\label{lem:convcrit}
Let $U\subset\R^N$ be open and convex, and let $v\in C(U)$. Suppose that $D^2\phi(x)\ge0$ whenever $x\in U$, $\phi$ is $C^2$ near $x$, $\phi\ge v$ near $x$ and $\phi(x)=v(x)$. Then $v$ is convex.
\end{lemma}

\begin{proof}
Suppose not. Then there are $a,b\in U$ such that, with $e=b-a$ and $L$ the affine function of $\scal{x-a}{e}$ with $L(a)=v(a)$, $L(b)=v(b)$,
\[
m:=\max_{0\le t\le1}\big(v-L\big)(a+te)\ >\ 0 .
\]
Write $x=a+te+x'$ with $x'\perp e$, and for $r>0$ let $C_r=\{x:\ 0\le t\le1,\ |x'|\le r\}$. Choose $r$ so small that $C_r\subset U$ and $|v-L|\le m/4$ on the two end discs $\{t=0\}$ and $\{t=1\}$ of $C_r$. Then choose $K$ so large that $Kr^2>\max_{C_r}|v-L|+m$. Consider
\[
W(x)=v(x)-L(x)+m\,(t-\tfrac12)^2-K|x'|^2\qquad\text{on }C_r .
\]
Then:
\begin{itemize}
\item on the segment $[a,b]$, $\max W\ge m$;
\item on the end discs, $W\le m/4+m/4<m$;
\item on the lateral boundary $|x'|=r$, $W<m$.
\end{itemize}
Hence $W$ attains its maximum over $C_r$ at an interior point $x_*$. The function
\[
\phi(x)=L(x)-m\,(t-\tfrac12)^2+K|x'|^2+W(x_*)
\]
satisfies $\phi\ge v$ near $x_*$ with equality at $x_*$, but $D^2\phi(x_*)[e,e]=-2m<0$. This is a contradiction.
\end{proof}

To analyse second derivatives of $H$ we remove the homogeneity in $w$. Since $H(\tau,\cdot)$ is positively $1$-homogeneous, its Hessian in $w$ is degenerate in the radial direction: \[
D^2_wH(\tau,w)\,w=0,\]
and $H$ is determined by its values on any hypersurface that meets each ray once. The unit sphere is the natural such hypersurface, but restricting $H$ to it produces first-order correction terms, as in $h_{\varphi\varphi}+h$ for $n=2$. We therefore restrict $H$ instead to an affine hyperplane, i.e.\ we use the central (gnomonic) projection. Fix a unit vector $\theta_0$ and put
\begin{align}
k(\tau,s)&=H(\tau,\theta_0+s),\qquad s\in\theta_0^\perp\cong\R^{n-1},\notag\\
H(\tau,w)&=\scal{w}{\theta_0}\;k\Big(\tau,\frac{w}{\scal{w}{\theta_0}}-\theta_0\Big)\qquad\text{for }\scal{w}{\theta_0}>0 .\label{eq:chart}
\end{align}
The second line is the inverse of the first: the affine hyperplane $\theta_0+\theta_0^\perp$ meets every ray of the half-space $\{\scal{w}{\theta_0}>0\}$ exactly once, and $H$ is recovered from $k$ by homogeneity. Since the hyperplane is affine, the restriction of a function to it has as Hessian the restriction of the Hessian: for $(t,e)\in\R\times\theta_0^\perp$,
\[
D^2k(\tau,s)[(t,e),(t,e)]=D^2H(\tau,\theta_0+s)[(t,e),(t,e)],
\]
with no first-order terms. Thus convexity properties of $H$ in $(\tau,w)$ near a ray correspond to convexity properties of $k$ in $(\tau,s)$, and the latter is a function on an open subset of $\R\times\R^{n-1}$ without any built-in degeneracy.

\begin{theorem}\label{thm:semiconvex}
Let $M$ be a compact connected $C^2$-immersed $p$-minimal hypersurface with $x(\partial M)\subset\Pi_{\tau_0}\cup\Pi_{\tau_1}$, and let $H$ be the support function \eqref{eq:Hdef} of the hulls of its sections. Let $\tau_0<a<b<\tau_1$, and put
\[
\varrho=\max_{x(M)}|y|,\qquad \Lambda=\frac{2\varrho}{\min\{a-\tau_0,\ \tau_1-b\}} .
\]
Then for every unit vector $\theta_0\in\R^n$ the function
\[
k(\tau,s)+\frac{\Lambda}{2}\big(\tau^2+|s|^2\big),
\]
where $k$ is given by \eqref{eq:chart}, is convex on $(a,b)\times\theta_0^\perp$. Consequently, $H$ is locally semiconvex on $(\tau_0,\tau_1)\times(\R^n\setminus0)$; in particular, $H$ is locally a difference of two convex functions.
\end{theorem}

\begin{proof}
We first obtain a uniform bound in the $\tau$-direction and then use the affine chart \eqref{eq:chart} to transfer upper test functions to $H$. Theorem~\ref{thm:main} then gives a uniform Hessian lower bound, from which semiconvexity follows after adding a quadratic correction.

\emph{Step 1: a bound for $H_\tau$.} We have $|H(\tau,\theta)|\le\varrho$ for $\tau\in[\tau_0,\tau_1]$ and $|\theta|=1$. Since $H(\cdot,\theta)$ is convex on $[\tau_0,\tau_1]$, for $\tau\in[a,b]$
\[
\partial^+_\tau H(\tau,\theta)\le\frac{H(\tau_1,\theta)-H(\tau,\theta)}{\tau_1-\tau}\le\Lambda,\qquad
\partial^-_\tau H(\tau,\theta)\ge\frac{H(\tau,\theta)-H(\tau_0,\theta)}{\tau-\tau_0}\ge-\Lambda .
\]

\emph{Step 2: transfer of test functions.} Let $\phi$ be $C^2$ near $(\tau_*,s_*)\in(a,b)\times\theta_0^\perp$, with $\phi\ge k$ near this point and equality at it. By \eqref{eq:chart},
\[
P(\tau,w)=\scal{w}{\theta_0}\;\phi\Big(\tau,\frac{w}{\scal{w}{\theta_0}}-\theta_0\Big)
\]
is $C^2$ and $1$-homogeneous, $P\ge H$ near $(\tau_*,\theta_0+s_*)$, with equality there, and $P(\tau,\theta_0+s)=\phi(\tau,s)$. By homogeneity, $P$ touches $H$ from above also at $(\tau_*,\theta)$, $\theta=(\theta_0+s_*)/\rho$, $\rho=|\theta_0+s_*|$. Since $P(\cdot,\theta)\ge H(\cdot,\theta)$ near $\tau_*$ with equality at $\tau_*$, and $H(\cdot,\theta)$ is convex,
\[
\partial^+_\tau H(\tau_*,\theta)\le P_\tau(\tau_*,\theta)\le\partial_\tau^-H(\tau_*,\theta)\le\partial^+_\tau H(\tau_*,\theta).
\]
Hence $H(\cdot,\theta)$ is differentiable at $\tau_*$, and by Step 1, $|P_\tau(\tau_*,\theta)|\le\Lambda$.

\emph{Step 3: the Hessian of $\phi$.} Let $(t,e)\in\R\times\theta_0^\perp$ and $\xi=e-\scal{e}{\theta}\theta\in\theta^\perp$. Since $P$ is $1$-homogeneous:
\begin{itemize}
\item $P_{\tau\tau}$ is $1$-homogeneous in $w$;
\item $\nabla_wP_\tau$ is $0$-homogeneous, with $\scal{\nabla_wP_\tau(\tau,\theta)}{\theta}=P_\tau(\tau,\theta)$ by Euler's identity;
\item $D^2_wP$ is $(-1)$-homogeneous, with $D^2_wP(\tau,\theta)\,\theta=0$.
\end{itemize}
Differentiating $\phi(\tau,s)=P(\tau,\theta_0+s)$ twice and evaluating at $(\tau_*,\theta)$, we get
\begin{align}
D^2\phi(\tau_*,s_*)[(t,e),(t,e)]
&=\rho\,P_{\tau\tau}\,t^2+2t\big(P_\tau\scal{\theta}{e}+\scal{b_P}{\xi}\big)+\rho^{-1}\scal{\mathcal R_P\,\xi}{\xi}\notag\\
&=Q_P\big(\sqrt\rho\,t,\ \xi/\sqrt\rho\big)+2t\,P_\tau\scal{\theta}{e},\label{eq:Hesschart}
\end{align}
where $b_P$ and $\mathcal R_P$ are as in Definition~\ref{def:visc}, and
\[
Q_P(t',\xi')=P_{\tau\tau}\,t'^2+2t'\scal{b_P}{\xi'}+\scal{\mathcal R_P\,\xi'}{\xi'} .
\]
By Theorem~\ref{thm:main}, $\mathcal R_P>0$ and $P_{\tau\tau}>\scal{\mathcal R_P^{-1}b_P}{b_P}$. So the Schur complement of $\mathcal R_P$ in the matrix of $Q_P$ is positive, and $Q_P$ is positive definite. Hence, by Step 2,
\[
D^2\phi(\tau_*,s_*)[(t,e),(t,e)]\ \ge\ -2\Lambda\,|t|\,|e|\ \ge\ -\Lambda\,(t^2+|e|^2).
\]

\emph{Step 4: convexity.} Let $v=k+\frac\Lambda2(\tau^2+|s|^2)$, which is continuous by Lemma~\ref{lem:convex}. If $\phi$ touches $v$ from above at a point of $(a,b)\times\theta_0^\perp$, then $\phi-\frac\Lambda2(\tau^2+|s|^2)$ touches $k$ from above there, so $D^2\phi\ge0$ by Step 3. Lemma~\ref{lem:convcrit} shows that $v$ is convex.

\emph{Step 5: semiconvexity of $H$.} We show that every point $(\bar\tau,\bar w)$ with $\bar\tau\in(\tau_0,\tau_1)$ and $\bar w\ne0$ has a convex neighbourhood $U$ and a constant $C$ such that $H+C(\tau^2+|w|^2)$ is convex on $U$.
Choose $\theta_0=\bar w/|\bar w|$ and $a<\bar\tau<b$ with $[a,b]\subset(\tau_0,\tau_1)$. Then $\scal{\bar w}{\theta_0}>0$. Let $U$ be a convex open neighbourhood of $(\bar\tau,\bar w)$ whose closure $\overline U$ is a compact subset of $(a,b)\times\{\scal{w}{\theta_0}>0\}$. On $\overline U$ put
\[
\gamma(w)=\scal{w}{\theta_0}>0,\qquad \Phi(\tau,w)=\Big(\tau,\frac{w}{\scal{w}{\theta_0}}-\theta_0\Big),\qquad q(\tau,s)=\frac\Lambda2(\tau^2+|s|^2).
\]
Both $\gamma$ and $\Phi$ are smooth on $\overline U$, and $\gamma$ is linear. By \eqref{eq:chart} and Step~4, $k=v-q$, so that
\begin{equation}\label{eq:Hsplit}
H=\gamma\cdot(v\circ\Phi)-\gamma\cdot(q\circ\Phi)\qquad\text{on }U .
\end{equation}

The function $\gamma\cdot(q\circ\Phi)$ is smooth on the compact set $\overline U$, so its Hessian is bounded there, say by $C_1$ in operator norm.
Next, let $K\subset(a,b)\times\theta_0^\perp$ be a compact convex set containing $\Phi(\overline U)$ in its interior. Being convex, $v$ is Lipschitz on $K$, say with constant $C_K$. It is also the upper envelope of its supporting affine functions at points of $K$: for $x\in\Phi(\overline U)$,
\[
v(x)=\max_{z\in K}\ell_z(x),\qquad \ell_z(x)=v(z)+\scal{p_z}{x-z},\quad p_z\in\partial v(z),
\]
with $|p_z|\le C_K$ for all $z\in K$. For each $z$, the function $f_z=\gamma\cdot(\ell_z\circ\Phi)$ is smooth on $\overline U$, and since $D^2\gamma=0$,
for every vector $X\in\R\times\R^n$,
\[
D^2f_z[X,X]=2\,\scal{\nabla\gamma}{X}\,\scal{p_z}{D\Phi\,X}+\gamma\,\scal{p_z}{D^2\Phi[X,X]} ,
\]
where $D^2\Phi[X,X]$ is the vector of second derivatives of the components of $\Phi$.
Hence $\|D^2f_z\|\le C_2$ on $\overline U$, with $C_2$ depending only on $C_K$ and on bounds for $\gamma$, $D\Phi$ and $D^2\Phi$ on $\overline U$, not on $z$. Consequently each $f_z+\frac{C_2}2(\tau^2+|w|^2)$ is convex on the convex set $U$. Since $\gamma>0$, we have \[\gamma\cdot(v\circ\Phi)=\sup_{z\in K}f_z \quad \text{on $U$},
\]
and a supremum of convex functions is convex. Thus $\gamma\cdot(v\circ\Phi)+\frac{C_2}2(\tau^2+|w|^2)$ is convex on $U$.

In summary, by \eqref{eq:Hsplit}, $H+C(\tau^2+|w|^2)$ is convex on $U$ with $C=\frac12(C_1+C_2)$. Hence $H$ is locally semiconvex, and on $U$ it is the difference of the two convex functions $H+C(\tau^2+|w|^2)$ and $C(\tau^2+|w|^2)$. For the calculus of semiconvex and semiconcave functions see \cite[\S2.1]{CannarsaSinestrari04}.

\end{proof}

\begin{remark}\label{rem:regularity}
(1) The semiconvexity is a consequence of the comparison with $M$, not of convexity alone. The proof uses only the qualitative part of \eqref{eq:mainP}, namely the positivity of $Q_P$. The constant $\Lambda$ enters only through the first-order term $2tP_\tau\scal{\theta}{e}$ in \eqref{eq:Hesschart}, which comes from the tilt of the chart; at $s_*=0$ it vanishes.

(2) The decomposition $k=\big(k+\frac\Lambda2(\tau^2+|s|^2)\big)-\frac\Lambda2(\tau^2+|s|^2)$ represents $H$, in every gnomonic chart, as a difference of convex functions. As far as we recall, a representation of this kind was the device used in \cite{T92}; here it is a consequence of the viscosity inequality rather than an assumption.
\end{remark}

\begin{corollary}\label{cor:pointwise}
Let $M$ and $H$ be as in Theorem~\ref{thm:semiconvex}. Then:
\begin{enumerate}[label=(\roman*)]
\item For almost every $(\tau,w)\in(\tau_0,\tau_1)\times(\R^n\setminus0)$, the function $H$ is twice differentiable at $(\tau,w)$ in the sense of Peano, $\mathcal R=D^2_wH(\tau,\theta)|_{\theta^\perp}$ is positive definite, $\theta=w/|w|$, and \eqref{eq:main} holds pointwise at $(\tau,\theta)$.
\item The distributional Hessian of $H$ in $(\tau,w)$ is a symmetric-matrix-valued Radon measure on $(\tau_0,\tau_1)\times(\R^n\setminus0)$. Its singular part is positive semidefinite, and the density of its absolutely continuous part coincides a.e.\ with the Peano Hessian in~(i).
\end{enumerate}
\end{corollary}

\begin{proof}
By Theorem~\ref{thm:semiconvex}, $H$ is locally semiconvex: near every point, $H+C(\tau^2+|w|^2)$ is convex for some $C$. For a convex function $F$, the distributional Hessian is a matrix of Radon measures with positive semidefinite values, and its absolutely continuous part is a.e.\ the Peano Hessian (Alexandrov's theorem) \cite[\S\S6.3--6.4]{EvansGariepy15}. Applying this to $F=H+C(\tau^2+|w|^2)$ gives (ii), since $D^2H$ and $D^2F$ differ by the constant matrix $2C\,\mathrm{Id}$, and the first assertion of (i).

At a point of (i), by homogeneity, $H$ is also twice differentiable at $(\tau,\theta)$. Let $T$ be its second-order Taylor polynomial there, and let $T_\varepsilon=T+\varepsilon(|\tau'-\tau|^2+|w'-\theta|^2)$. Then $T_\varepsilon$ touches $H$ from above near $(\tau,\theta)$. So does its $1$-homogeneous extension $P(\tau',w')=|w'|\,T_\varepsilon(\tau',w'/|w'|)$. At $(\tau,\theta)$ the quantities $P_\tau$, $P_{\tau\tau}$, $b_P$ and $\mathcal R_P$ of Definition~\ref{def:visc} are those of $T_\varepsilon$. For $\mathcal R_P$ this uses $\partial_rT_\varepsilon=\scal{\nabla_wH}{\theta}=H=T_\varepsilon$ at $(\tau,\theta)$, by Euler's identity. Definition~\ref{def:visc} applied to $P$ gives $\mathcal R+2\varepsilon\,\mathrm{Id}>0$ and \eqref{eq:mainP} with $P_{\tau\tau}=H_{\tau\tau}+2\varepsilon$, $\mathcal R_P=\mathcal R+2\varepsilon\,\mathrm{Id}$ and $b_P=b$. If $\mathcal R$ had a zero eigenvalue, then $\tr(\mathcal R+2\varepsilon\,\mathrm{Id})^{-1}\to\infty$ as $\varepsilon\to0$, while the left side stays bounded. Hence $\mathcal R>0$, and letting $\varepsilon\to0$ gives \eqref{eq:main}.
\end{proof}

\section{Surfaces in \texorpdfstring{$\R^3$}{R3}: a Monge--Amp\`ere inequality}\label{sec:MA}

Throughout this section we assume $n=2$ and write
\[
h(\tau,\varphi)=H(\tau,\theta(\varphi)), \qquad \theta(\varphi)=(\cos\varphi,\sin\varphi).
\]

\begin{theorem}\label{thm:MA}
$h$ is a viscosity subsolution of \eqref{eq:MA} on $(\tau_0,\tau_1)\times\R$. That is, if $q\in C^2$ near $(\tau,\varphi)$, $q\ge h$ near $(\tau,\varphi)$ and $q(\tau,\varphi)=h(\tau,\varphi)$, then $q_{\varphi\varphi}+q>0$ and
\[
q_{\tau\tau}(q_{\varphi\varphi}+q)-q_{\tau\varphi}^2\ \ge\ \frac{1+q_\tau^2}{p-1}\quad\text{at }(\tau,\varphi).
\]
\end{theorem}

\begin{proof}
Put $P(\tau',w)=|w|\,q(\tau',\arg w)$. Then $P$ is $C^2$, $1$-homogeneous and touches $H$ from above at $(\tau,\theta(\varphi))$. With $u=\theta'(\varphi)$ we have $\mathcal R_P=q_{\varphi\varphi}+q$ by \ref{S4}, $b_P=\partial_u P_\tau=q_{\tau\varphi}$, $P_\tau=q_\tau$ and $P_{\tau\tau}=q_{\tau\tau}$. Multiplying \eqref{eq:mainP} by $\mathcal R_P>0$ gives the claim.
\end{proof}

The left side of \eqref{eq:MA} is the determinant of the matrix
\(
\begin{pmatrix} h_{\tau\tau}&h_{\tau\varphi}\\ h_{\tau\varphi}&h_{\varphi\varphi}+h\end{pmatrix},
\)
the Hessian of the $1$-homogeneous extension $H$ restricted to $\mathrm{span}\{\partial_\tau,\theta'\}$. Inequality \eqref{eq:MA} is therefore of Monge--Amp\`ere type with a gradient-dependent right-hand side; see \cite{Gutierrez01} for the theory of such equations.

By Theorem~\ref{thm:semiconvex}, Theorem~\ref{thm:MA} also holds in the sense of measures, and the singular part of this matrix records the edges of the hulls.

\begin{proposition}[Theorem~C in the sense of measures]\label{prop:MAmeasure}
Let $n=2$. The distribution
\[
\mathbb M=\begin{pmatrix} h_{\tau\tau}&h_{\tau\varphi}\\ h_{\tau\varphi}&h_{\varphi\varphi}+h\end{pmatrix}
\]
on $(\tau_0,\tau_1)\times\R$ is a symmetric-matrix-valued Radon measure. Write $\mathbb M=\mathbb M^{ac}\,d\tau\,d\varphi+\mathbb M^{s}$ for its Lebesgue decomposition. Then:
\begin{enumerate}[label=(\roman*)]
\item $\mathbb M^s$ is positive semidefinite;
\item for a.e.\ $(\tau,\varphi)$ the matrix $\mathbb M^{ac}(\tau,\varphi)$ is positive definite and
\[
\det\mathbb M^{ac}\ \ge\ \frac{1+h_\tau^2}{p-1};
\]
\item the entry $\mathbb M_{22}$ disintegrates into the area functions $S_\tau=S(\Omega(\tau),\cdot)$ of \ref{S6}:
\[
\int\zeta\,d\mathbb M_{22}=\int_{\tau_0}^{\tau_1}\Big(\int_{\R}\zeta(\tau,\cdot)\,dS_\tau\Big)d\tau\qquad\text{for }\zeta\in C_c\big((\tau_0,\tau_1)\times\R\big),
\]
where $S_\tau$ is regarded as a $2\pi$-periodic measure on $\R$.
\end{enumerate}
In particular $\mathbb M\ge0$, and
\[
\mathbb M^s_{22}\ \ge\ d\tau\otimes\sum_{E}|E|\,\delta_{\theta_E},
\]
where, for each $\tau$, the sum runs over the edges $E$ of $\Omega(\tau)$, with outer normals $\theta_E$.
\end{proposition}

\begin{proof}
The map $(\tau,\varphi)\mapsto(\tau,\theta(\varphi))$ is smooth, so $h$ is locally semiconvex by Theorem~\ref{thm:semiconvex}, by the argument of its Step~5. Hence the distributional Hessian $D^2h$ in $(\tau,\varphi)$ is a matrix-valued Radon measure with positive semidefinite singular part. Since $\mathbb M=D^2h+\mathrm{diag}(0,h)$ and $h$ is continuous, $\mathbb M^s=(D^2h)^s\ge0$, which is (i).

By Alexandrov's theorem, at a.e.\ $(\tau,\varphi)$ the function $h$ has a second-order Peano expansion with Hessian $\mathbb M^{ac}-\mathrm{diag}(0,h)$. Let $T$ be the corresponding Taylor polynomial and $q_\varepsilon=T+\varepsilon(|\tau'-\tau|^2+|\varphi'-\varphi|^2)$, which touches $h$ from above near $(\tau,\varphi)$. Theorem~\ref{thm:MA} applied to $q_\varepsilon$, followed by $\varepsilon\to0$, gives $\mathbb M^{ac}_{22}\ge0$ and $\det\mathbb M^{ac}\ge(1+h_\tau^2)/(p-1)>0$. Hence $\mathbb M^{ac}$ is positive definite, which is (ii).

For (iii), let first $\zeta\in C^2_c$. By Fubini,
\[
\int h\,(\zeta_{\varphi\varphi}+\zeta)\,d\tau\,d\varphi=\int_{\tau_0}^{\tau_1}\Big(\int h(\tau,\cdot)\big(\zeta_{\varphi\varphi}(\tau,\cdot)+\zeta(\tau,\cdot)\big)d\varphi\Big)d\tau ,
\]
and the inner integral is $\int\zeta(\tau,\cdot)\,dS_\tau$ by \ref{S6}. Both sides of (iii) are positive linear functionals of $\zeta$, and $\tau\mapsto S_\tau(\R/2\pi\mathbb Z)$ is locally bounded, so (iii) extends to $\zeta\in C_c$ by density. The last assertion follows from (iii): the measure $d\tau\otimes\sum_E|E|\delta_{\theta_E}$ is the atomic part of the disintegration, and it is singular with respect to $d\tau\,d\varphi$.
\end{proof}

Thus the edges of the hulls live in the singular part of $\mathbb M$, which carries no information from \eqref{eq:MA}, while \eqref{eq:MA} governs the absolutely continuous part. By Lemma~\ref{lem:smooth}, $S_\tau$ charges every open arc. In the setting of Section~\ref{sec:MAC}, the edges are unions of chords, so $\mathbb M_{22}^s$ dominates the total chord length (Corollary~\ref{cor:C1}).

\begin{proposition}[Equality case]\label{prop:equality}
Let $M\subset\R^3$ be a $C^3$ $p$-minimal surface such that $x_3$ has no critical points on $M^\circ$ and every intermediate section $\Sigma(\tau)$ is a simple closed curve of positive curvature, so that $\Sigma(\tau)=\partial\Omega(\tau)$. Then $h\in C^2$ and equality holds in \eqref{eq:MA} on $(\tau_0,\tau_1)\times\R$.
\end{proposition}

\begin{proof}
By assumption $\Sigma(\tau)=\partial\Omega(\tau)$, and the point $y(\tau,\theta)$ of $\Sigma(\tau)$ with outer normal $\theta$ depends in a $C^2$ way on $(\tau,\theta)$ by the implicit function theorem. Hence $H(\tau,w)=\scal{y(\tau,w/|w|)}{w}$ is $C^2$, and \eqref{eq:main} holds pointwise by Theorem~\ref{thm:main}. Let $m$ be a contact point and choose $w(m')$ to be the outer unit normal of $\Sigma(\tau(m'))$ at $y(m')$. Then $\psi\equiv0$ near $m$, so equality holds in \eqref{eq:Phi} for this particular choice of $W$. Since \eqref{eq:mainP} is equivalent to $\sup_W(\text{right side of }\eqref{eq:Phi})\le0$ and this supremum is attained, the supremum equals $0$, which is equality in \eqref{eq:main}.
\end{proof}

\begin{example}[Catenoid]
For $x_1^2+x_2^2=\cosh^2x_3$ we have $h=\cosh\tau$, and both sides of \eqref{eq:MA} with $p=2$ equal $\cosh^2\tau$.
\end{example}

\begin{example}[Circular sections]\label{ex:circles}
Let $h=\scal{\xi(\tau)}{\theta}+R(\tau)$, i.e. $\Omega(\tau)$ is the disc of radius $R$ centred at $\xi$. Then $h_{\varphi\varphi}+h=R$, $h_\tau=R'+\scal{\xi'}{\theta}$ and $h_{\tau\varphi}=\scal{\xi'}{\theta'}$. Using $\scal{\xi'}{\theta}^2+\scal{\xi'}{\theta'}^2=|\xi'|^2$, equality in \eqref{eq:MA} with $p=2$ for all $\varphi$ becomes
\[
RR''-1-R'^2-|\xi'|^2+\scal{\theta}{R\xi''-2R'\xi'}=0\qquad\forall\theta,
\]
that is,
\begin{equation}\label{eq:Riemann}
RR''=1+R'^2+|\xi'|^2,\qquad R\xi''=2R'\xi' \ \ (\text{equivalently } \xi'=cR^2,\ c\in\R^2).
\end{equation}
This is the system describing the catenoid ($c=0$) and Riemann's examples ($c\ne0$), the only minimal surfaces foliated by circles in parallel planes \cite{Riemann1867,Enneper1869}. For $p\ne2$, equality for all $\theta$ forces $\xi'=0$ and $RR''=(1+R'^2)/(p-1)$.
\end{example}

Thus for two circles $\Gamma_0,\Gamma_1$ the situation is as follows. By Shiffman \cite{Shiffman56}, every minimal annulus spanning them has circular sections and therefore solves \eqref{eq:Riemann}. By Theorem~\ref{thm:MA}, the hull data $h$ of \emph{any} connected minimal surface spanning them, of any genus, is a viscosity subsolution of the same equation with the same boundary values $h(\tau_i,\cdot)$.

\section{Integrated consequences}\label{sec:integrated}

In this section we average Theorem~\ref{thm:main} over the sphere of directions. The result is a closed ordinary differential inequality for the mean half-width of the hulls of the sections. Its equality case is the catenoid for $n=p=2$, and it contains several known nonexistence and lifetime results in a sharper form.

Integration requires the inequality \eqref{eq:main} pointwise almost everywhere. By Corollary~\ref{cor:pointwise}, this holds without any assumption, as a consequence of the semiconvexity of $H$ (Theorem~\ref{thm:semiconvex}). For each $\tau$ let $\sigma_\tau=\sigma(\Omega(\tau),\cdot)$ be the measure of \ref{S6}, and let $\sigma_\tau=\sigma_\tau^{ac}+\sigma_\tau^{s}$ be its Lebesgue decomposition with respect to $d\theta$.

\subsection{The mean half-width}
Recall the \emph{mean half-width} \eqref{eq:Bdef} of the hull of $\Sigma(\tau)$:
\[
B(\tau)=\frac{1}{|S^{n-1}|}\int_{S^{n-1}}H(\tau,\theta)\,d\theta .
\]
It is translation invariant, and $B(\tau)\le R(\tau)$, the circumradius of $\Sigma(\tau)$, by \ref{S2}. For $n=2$, $2\pi B(\tau)=L(\tau)$ is the perimeter of $\Omega(\tau)$ by Cauchy's formula \cite[\S1.7]{Schneider14}. By Lemma~\ref{lem:convex}, $B$ is convex on $[\tau_0,\tau_1]$.
For $n\ge2$ put
\[
\lambda(\tau)=\frac{\sigma_\tau^s(S^{n-1})}{(n-1)\,|S^{n-1}|},
\]
the normalized mass of the singular part of $\sigma_\tau$. By \ref{S6}, $0\le\lambda(\tau)\le B(\tau)$. For $n=2$, $\sigma_\tau=S_\tau$ is the area function, so $2\pi\lambda(\tau)$ is at least the total length $\ell(\tau)$ of the edges of $\Omega(\tau)$.

\begin{proposition}\label{prop:Bineq}
Let $\beta=\frac{n-1}{p-1}$. Then $\lambda$ is measurable, $\lambda<B$ a.e., and, in the sense of measures on $(\tau_0,\tau_1)$,
\begin{equation}\label{eq:Bineqlambda}
B''\ \ge\ \beta\,\frac{1+B'^2}{B-\lambda}.
\end{equation}
In particular,
\begin{equation}\label{eq:Bineq}
B''\ \ge\ \beta\,\frac{1+B'^2}{B}.
\end{equation}
\end{proposition}

\begin{proof}
Let $h=H|_{(\tau_0,\tau_1)\times S^{n-1}}$. By Corollary~\ref{cor:pointwise}(i), at a.e.\ $(\tau,\theta)$ the matrix $\mathcal R=D^2_wH|_{\theta^\perp}$ is positive definite and \eqref{eq:main} holds. By the arithmetic--harmonic mean inequality,
\[
\tr\mathcal R^{-1}\ge(n-1)^2/\tr\mathcal R.
\]
Dropping the nonnegative last term of \eqref{eq:main},
\begin{equation}\label{eq:pointwiseB}
h_{\tau\tau}\ \ge\ \frac{(n-1)^2}{p-1}\,\frac{1+h_\tau^2}{\tr\mathcal R}\qquad\text{a.e. in }(\tau_0,\tau_1)\times S^{n-1}.
\end{equation}
We integrate over $S^{n-1}$, using three facts:
\begin{itemize}
\item By Cauchy--Schwarz,
\[
\int_{S^{n-1}}\frac{1+h_\tau^2}{\tr\mathcal R}\,d\theta\ \ge\ \frac{\big(\int\sqrt{1+h_\tau^2}\,d\theta\big)^2}{\int\tr\mathcal R\,d\theta}.
\]
\item By Minkowski's integral inequality,
\[
\int_{S^{n-1}}\sqrt{1+h_\tau^2}\,d\theta\ \ge\ \Big(|S^{n-1}|^2+\Big(\int h_\tau\,d\theta\Big)^2\Big)^{1/2}=|S^{n-1}|(1+B'^2)^{1/2}.
\]
Here we used $B'(\tau)=|S^{n-1}|^{-1}\int h_\tau(\tau,\theta)\,d\theta$ for a.e.\ $\tau$, which follows by dominated convergence from the uniform bound $|\partial^\pm_\tau H|\le\Lambda$ of Step~1 in the proof of Theorem~\ref{thm:semiconvex}.
\item For a.e.\ $\tau$, $H$ is twice differentiable at a.e.\ $(\tau,\theta)$ by Fubini's theorem, so $\tr\mathcal R(\tau,\cdot)$ is a.e.\ the trace of the Alexandrov Hessian of the convex function $H(\tau,\cdot)$ on $\theta^\perp$. It is therefore the density of $\sigma^{ac}_\tau$ \cite[\S6.4]{EvansGariepy15}, and by \ref{S6}
\[
\int_{S^{n-1}}\tr\mathcal R\,d\theta=\sigma_\tau(S^{n-1})-\sigma_\tau^s(S^{n-1})=(n-1)\,|S^{n-1}|\,\big(B-\lambda\big).
\]
\end{itemize}
The left side of the last identity is positive since $\mathcal R>0$ a.e., so $\lambda<B$ for a.e.\ $\tau$; it is measurable in $\tau$ by Fubini, hence so is $\lambda$. Together these give, for a.e.\ $\tau$,
\[
\frac{1}{|S^{n-1}|}\int_{S^{n-1}}h_{\tau\tau}\,d\theta\ \ge\ \beta\,\frac{1+B'^2}{B-\lambda}.
\]
Finally, for each $\theta$ the function $h(\cdot,\theta)$ is convex, so its distributional second derivative dominates its pointwise second derivative. Integrating over $\theta$ by Fubini, $B''$ dominates, as a measure, the left-hand side above.
\end{proof}

Thus the singular part of the area measures, in particular the edges of the hulls for $n=2$, makes the inequality strictly stronger. We return to this in Remark~\ref{rem:strict}.

\begin{remark}
The averaging has a conceptual explanation for $n=2$. By Theorem~\ref{thm:MA}, $\sqrt{\det}$ of the matrix $\begin{psmallmatrix}h_{\tau\tau}&h_{\tau\varphi}\\h_{\tau\varphi}&h_{\varphi\varphi}+h\end{psmallmatrix}$ is at least $\sqrt{(1+h_\tau^2)/(p-1)}$. The left side is concave in the matrix and the right side is convex in $h_\tau$, so the inequality survives averaging over rotations. The rotation average of $\Omega(\tau)$ is the disc of radius $B(\tau)$, for which the matrix is $\mathrm{diag}(B'',B)$.
\end{remark}

\subsection{Comparison with generalized catenaries}
For $\beta>0$ let 
\[
G_\beta(s)=\int_1^s\frac{dx}{\sqrt{x^{2\beta}-1}}, \qquad s\ge1. 
\]
Note that $G_1=\operatorname{arccosh}$, and $\kappa_\beta=G_\beta(\infty)<\infty$ if and only if $\beta>1$. The equality case of \eqref{eq:Bineq} has the first integral 
\[
1+B'^2=(B/m)^{2\beta}, \qquad m>0. 
\]
Its solutions are the generalized catenaries $B=m\,C_\beta((\tau-c)/m)$, where $C_\beta$ is the even solution with $C_\beta(0)=1$, inverse to $G_\beta$; it is defined on $[0,\infty)$ if $\beta\le1$ and on $[0,\kappa_\beta)$ if $\beta>1$. For $\beta=1$, $C_1=\cosh$, and $2\pi B=L$ is the perimeter of the sections of a catenoid.

\begin{lemma}\label{lem:ODE}
Let $B>0$ be convex on $(a,b)$ and satisfy \eqref{eq:Bineq} in the sense of measures. Then $B$ is strictly convex, and its minimum $B_{\min}$ on $[a,b]$ is attained at a unique point $\tau_m$, with $B$ extended continuously to $[a,b]$. Moreover,
\begin{equation}\label{eq:catenary}
|\tau-\tau_m|\ \le\ B_{\min}\,G_\beta\big(B(\tau)/B_{\min}\big),\qquad\text{equivalently}\qquad B(\tau)\ \ge\ B_{\min}\,C_\beta\Big(\frac{\tau-\tau_m}{B_{\min}}\Big).
\end{equation}
\end{lemma}

\begin{proof}
By \eqref{eq:Bineq}, $B''\ge\beta/B>0$, so $B$ is strictly convex and $\tau_m$ is unique. We treat $\tau>\tau_m$; the other side is symmetric. There, $B'\ge0$, $B'$ is nondecreasing, and $B$ is increasing. Let $\Phi(s)=(1+s^2)^{-1/(2\beta)}$, which is nonincreasing on $[0,\infty)$, and put $\mu=B\,\Phi(B')$.

Decompose the measure $dB'$ into its absolutely continuous, jump and Cantor parts, and use the chain rule for functions of bounded variation. On the absolutely continuous part, by \eqref{eq:Bineq},
\[
d\mu=B'\Phi(B')\,d\tau+B\,\Phi'(B')\,dB'\le\Big(B'\Phi(B')-\frac{B'}{1+B'^2}\,\Phi(B')\,(1+B'^2)\Big)d\tau=0,
\]
using $\Phi'(s)=-\frac{s}{\beta}\frac{\Phi(s)}{1+s^2}$. At a jump of $B'$ upwards, $\Phi(B')$ jumps down. The singular continuous part of $dB'$ is nonnegative and $\Phi'\le0$. Hence $\mu$ is nonincreasing on $[\tau_m,b)$, and 
\[
\mu\le B_{\min}\Phi(B'(\tau_m+))\le B_{\min}.
\]

Thus $1+B'^2\ge(B/B_{\min})^{2\beta}$, i.e. $B'\ge\sqrt{(B/B_{\min})^{2\beta}-1}$ a.e. on $(\tau_m,b)$. Since $B$ is increasing and absolutely continuous, the change of variables $x=B/B_{\min}$ gives $\tau-\tau_m\le B_{\min}G_\beta(B(\tau)/B_{\min})$, which is the claim.
\end{proof}

\begin{theorem}\label{thm:integrated}
Let $B_i=B(\tau_i)$ and $d=\tau_1-\tau_0$. Then $B(\tau)\ge B_{\min}C_\beta((\tau-\tau_m)/B_{\min})$ for all $\tau$, and:
\begin{enumerate}[label=(\roman*)]
\item if $\beta>1$, i.e. $n>p$, then $d\le2\kappa_\beta B_{\min}\le2\kappa_\beta\min_\tau R(\tau)$;
\item if $\beta=1$, then
\[
d\ \le\ D(B_0,B_1):=\max_{0<m\le\min\{B_0,B_1\}}m\Big[\operatorname{arccosh}\frac{B_0}{m}+\operatorname{arccosh}\frac{B_1}{m}\Big];
\]
\item if $\beta<1$, i.e. $p>n$, then $B$ grows at least like $C_\beta$, and $C_\beta(t)\sim c\,t^{1/(1-\beta)}$ as $t\to\infty$.
\end{enumerate}
\end{theorem}

\begin{proof}
Apply Lemma~\ref{lem:ODE} with Proposition~\ref{prop:Bineq}. For (i) and (ii), $d=(\tau_1-\tau_m)+(\tau_m-\tau_0)$; bound each term by \eqref{eq:catenary}, and use $G_\beta\le\kappa_\beta$ for (i). For (ii), $B_{\min}\le\min\{B_0,B_1\}$, which gives the maximum over $m$. For (iii), the first integral gives $C_\beta'\sim C_\beta^{\beta}$, whence the growth rate.
\end{proof}

The number $D(B_0,B_1)$ is the largest height of a piece of a catenoid of revolution bounded by two coaxial circles of radii $B_0$ and $B_1$. In particular $D(R,R)=c_*R$, where
\begin{equation}\label{xnoll}
c_*=\max_{0<s\le1}2s\,\operatorname{arccosh}(1/s)=\frac{2}{\sinh x_0}\approx1.3255
\end{equation}
is the classical critical ratio for catenoids. Here $x_0\approx1.1997$ is the positive root of $x\tanh x=1$: with $s=1/\cosh x$, the function becomes $2x/\cosh x$, whose critical point satisfies $x\tanh x=1$, and then $\cosh x_0=x_0/\sqrt{x_0^2-1}$, so that $2x_0/\cosh x_0=2\sqrt{x_0^2-1}=2/\sinh x_0$. The extremal piece, between two circles of radius $R$ at height $c_*R$, is the marginally stable catenoid.

\subsection{Comparison with known results}

\begin{corollary}\label{cor:BB}
Let $M\subset\R^3$ be a compact connected immersed minimal surface with $\partial M=\Gamma_0\cup\Gamma_1$, $\Gamma_i\subset\Pi_{\tau_i}$ closed curves. Let $\rho_i=\operatorname{Per}(\conv\Gamma_i)/(2\pi)$. Then
\[
\tau_1-\tau_0\ \le\ D(\rho_0,\rho_1).
\]
Equivalently, if $\tau_1-\tau_0>D(\rho_0,\rho_1)$, then $\Gamma_0\cup\Gamma_1$ bounds no connected immersed minimal surface, of any genus.
\end{corollary}

\begin{proof}
Here $n=p=2$, so $\beta=1$. By the convex hull property \cite[Cor.~3]{T97}, $\Omega(\tau_i)=\conv\Gamma_i$, so $B_i=\rho_i$. Apply Theorem~\ref{thm:integrated}(ii).
\end{proof}

The equality case of Corollary~\ref{cor:BB} is rigid.

\begin{theorem}[Rigidity]\label{thm:rigidity}
In Corollary~\ref{cor:BB}, assume in addition that $x$ is an immersion up to the boundary and that $x|_{\partial M}$ is injective, i.e. each $\Gamma_i$ is a finite union of pairwise disjoint Jordan curves. If
\[
\tau_1-\tau_0=D(\rho_0,\rho_1),
\]
then $\Gamma_0$ and $\Gamma_1$ are coaxial circles of radii $\rho_0$ and $\rho_1$, and $x$ maps $M$ diffeomorphically onto the piece of the catenoid of revolution
\[
K=\Big\{(y,\tau):\ |y-\xi_0|=m_*\cosh\frac{\tau-\tau_m}{m_*},\ \tau_0\le\tau\le\tau_1\Big\},
\]
where $m_*$ is a maximizer in the definition of $D(\rho_0,\rho_1)$, $\xi_0\in\R^2$ and $\tau_m\in(\tau_0,\tau_1)$. In particular, $K$ contains the neck of the catenoid.
\end{theorem}

\begin{proof}
Here $n=p=2$, $\beta=1$ and $B_i=\rho_i$. Write $d=\tau_1-\tau_0$.

\emph{Step 1: the mean half-width is a catenary.} Let $\tilde B$ be the continuous extension of $B|_{(\tau_0,\tau_1)}$ to $[\tau_0,\tau_1]$; by convexity $\tilde B(\tau_i)\le B(\tau_i)=B_i$. Lemma~\ref{lem:ODE} applies to $\tilde B$, and the chain below with $\tilde B(\tau_i)$ in place of $B_i$ ends with $D(\tilde B(\tau_0),\tilde B(\tau_1))\le D(B_0,B_1)$. Since $D$ is strictly increasing in each argument, equality forces $\tilde B(\tau_i)=B_i$. Thus $B$ is convex and continuous on $[\tau_0,\tau_1]$, hence absolutely continuous there. In the proof of Theorem~\ref{thm:integrated}(ii),
\[
d=(\tau_1-\tau_m)+(\tau_m-\tau_0)\ \le\ B_{\min}\Big[\operatorname{arccosh}\frac{B_1}{B_{\min}}+\operatorname{arccosh}\frac{B_0}{B_{\min}}\Big]\ \le\ D(B_0,B_1),
\]
so both inequalities are equalities, and $m_*=B_{\min}$ is a maximizer. Since 
\[
\frac{d}{dm}\big[m\operatorname{arccosh}(\rho/m)\big]=
-\frac{\rho}{m\sqrt{\rho^2-m^2}}
\to-\infty
\] 
as $m\uparrow\rho$, a maximizer satisfies $m_*<\min\{\rho_0,\rho_1\}$, so $B_{\min}<B_i$ and $\tau_m\in(\tau_0,\tau_1)$. On $(\tau_m,\tau_1)$ we have $B'>0$ and, by Lemma~\ref{lem:ODE}, $B'\ge\sqrt{(B/B_{\min})^2-1}$ a.e. Hence for $\tau\in(\tau_m,\tau_1]$
\[
\tau-\tau_m=\int_{B_{\min}}^{B(\tau)}\frac{dB}{B'}\ \le\ B_{\min}\operatorname{arccosh}\frac{B(\tau)}{B_{\min}},
\]
with equality at $\tau=\tau_1$. Equality at $\tau_1$ forces $B'=\sqrt{(B/B_{\min})^2-1}$ a.e.\ on $(\tau_m,\tau_1)$, and then equality holds for every $\tau$. Arguing in the same way on $(\tau_0,\tau_m)$,
\[
B(\tau)=m_*\cosh\frac{\tau-\tau_m}{m_*},\qquad\tau\in[\tau_0,\tau_1].
\]
In particular $B$ is smooth and $B''=(1+B'^2)/B$.

\emph{Step 2: the hulls are coaxial discs.} The proof of Proposition~\ref{prop:Bineq} gives, as measures,
\[
B''\ \ge\ A(\tau)\,d\tau,\qquad A(\tau)=\frac1{2\pi}\int_0^{2\pi}h_{\tau\tau}\,d\varphi ,
\]
with $h_{\tau\tau}$ the pointwise second derivative. For a.e.\ $\tau$ the proof of Proposition~\ref{prop:Bineq} gives the chain
\[
A\ \overset{(a)}{\ge}\ \frac1{2\pi}\int\frac{1+h_\tau^2}{r}\,d\varphi\ \overset{(b)}{\ge}\ \frac1{2\pi}\,\frac{\big(\int\sqrt{1+h_\tau^2}\,d\varphi\big)^2}{\int r\,d\varphi}\ \overset{(c)}{\ge}\ \frac{1+B'^2}{B-\lambda}\ \overset{(d)}{\ge}\ \frac{1+B'^2}{B},
\]
where $r=h_{\varphi\varphi}+h$ is the pointwise radius of curvature. Here:
\begin{itemize}
\item (a) is the pointwise inequality \eqref{eq:MA} with $p=2$, $h_{\tau\tau}\ge(1+h_\tau^2)/r+h_{\tau\varphi}^2/r$, integrated in $\varphi$ after dropping the last term;
\item (b) is the Cauchy--Schwarz inequality;
\item (c) is Minkowski's inequality for the vectors $(1,h_\tau)$, combined with the identity $\int r\,d\varphi=2\pi(B-\lambda)$;
\item (d) uses $\lambda\ge0$.
\end{itemize}
We claim that (a)--(d) are equalities for a.e.\ $\tau$. By Step~1, $B$ is smooth with $B''=(1+B'^2)/B$, while the proof of Proposition~\ref{prop:Bineq} gives $B''\ge A\,d\tau$ as measures. Hence $A\le(1+B'^2)/B$ for a.e.\ $\tau$. Together with the chain this gives $A=(1+B'^2)/B$, so all intermediate inequalities are equalities. 
Fix such a $\tau$.
\begin{itemize}
\item Equality in (d) gives $\lambda(\tau)=0$, so $S_\tau=r(\tau,\cdot)\,d\varphi$ has no singular part.
\item Equality in (c) means equality in Minkowski's inequality, i.e. the vectors $(1,h_\tau(\tau,\varphi))$ are parallel for a.e.\ $\varphi$. Since their first components are equal, $h_\tau(\tau,\cdot)$ is a.e.\ constant.
\item Equality in (b) means that $\sqrt{1+h_\tau^2}/r$ is a.e.\ independent of $\varphi$. Since the numerator is, so is $r$.
\end{itemize}
Thus $S_\tau=c\,d\varphi$ for a constant $c$, i.e. $h_{\varphi\varphi}+h=c$ in the sense of distributions, and $c=B(\tau)$ because the total mass of $S_\tau$ is $2\pi B(\tau)$. Hence $h(\tau,\varphi)=B(\tau)+\scal{\xi(\tau)}{\theta(\varphi)}$ for some $\xi(\tau)\in\R^2$, i.e. $\Omega(\tau)$ is the disc of radius $B(\tau)$ centred at $\xi(\tau)$.

This holds for a.e.\ $\tau$, hence for all $\tau\in(\tau_0,\tau_1)$ by continuity of $h$ and $B$; in particular $\xi$ is continuous there. At an endpoint, say $\tau_0$, the discs $\Omega(\tau)$ converge as $\tau\downarrow\tau_0$ to a disc with support function $g$. By convexity in $\tau$, $h(\tau_0,\cdot)\ge g$, and both have mean $B_0$ by Step~1. Hence $h(\tau_0,\cdot)=g$ a.e., and everywhere by continuity in $\varphi$. So $\Omega(\tau_0)$ is a disc too, and $\xi$ extends continuously to $[\tau_0,\tau_1]$. Each component of $\xi$ is a difference of two functions convex in $\tau$, e.g. $\xi_1=h(\cdot,0)-B$, so $\xi$ is locally Lipschitz and differentiable a.e. At a.e.\ $\tau$ where $\xi$ is differentiable,
\[
h_\tau(\tau,\varphi)=B'(\tau)+\scal{\xi'(\tau)}{\theta(\varphi)}\quad\text{for all }\varphi .
\]
The left side is a.e.\ independent of $\varphi$, so $\xi'(\tau)=0$. Hence $\xi\equiv\xi_0$ is constant.

\emph{Step 3: $x(M)$ contains a catenoid.} By Steps~1 and~2 we conclude that $\Omega(\tau)$ is the disc of radius $m_*\cosh((\tau-\tau_m)/m_*)$ centred at $\xi_0$, for every $\tau\in[\tau_0,\tau_1]$. Every boundary point of a disc is an extreme point, so $\partial\Omega(\tau)\subset\Sigma(\tau)$ by \ref{S5}. Hence $K\subset x(M)$. At the end levels we use $\Sigma(\tau_i)=\Gamma_i$, which holds by the strong maximum principle $x_3(M^\circ)\subset(\tau_0,\tau_1)$. At the ends, $\Gamma_i\subset\Omega(\tau_i)$ and $\Gamma_i\supset\partial\Omega(\tau_i)=:C_i$.

\emph{Step 4: analytic continuation.} Let $F(y,t)=|y-\xi_0|^2-m_*^2\cosh^2((t-\tau_m)/m_*)$. It is real-analytic on $\R^3$, and $K=\{F=0,\ \tau_0\le t\le\tau_1\}$. It is well known that the immersion $x$ is real-analytic on $M^\circ$ in conformal coordinates \cite{Osserman86,DHS10}, so $f=F\circ x$ is real-analytic on the connected surface $M^\circ$.

Suppose $f\not\equiv0$. Then its zero set $Z$ has measure zero in every coordinate chart. Since $x$ is locally Lipschitz, $x(Z)$ has zero two-dimensional Hausdorff measure. But $K\cap\{\tau_0<t<\tau_1\}\subset x(Z)$, because $x_3(\partial M)\subset\{\tau_0,\tau_1\}$, and this set has positive area. This is a contradiction, so $f\equiv0$ and $x(M)\subset K$.

In summary we have: by the strong maximum principle, $x$ maps $M^\circ$ into the open annulus $K^\circ=K\cap\{\tau_0<t<\tau_1\}$ and $\partial M$ into $C_0\cup C_1$. So $\Gamma_i\subset C_i$, hence $\Gamma_i=C_i$. The map $x\colon M^\circ\to K^\circ$ is an immersion between surfaces, i.e. a local diffeomorphism. 
It is proper, because the preimage of a compact subset of $K^\circ$ is closed in $M$ and disjoint from $\partial M$. A proper local diffeomorphism between connected manifolds is a finite covering map, see e.g.\ \cite[Prop.~4.46]{Lee13}; hence $x\colon M^\circ\to K^\circ$ is a covering of some degree $k\ge1$. Near a point of $\partial M$, the immersion $x$ is injective on a half-disc and maps it onto a half-disc of $K$ along $C_0\cup C_1$. Hence $x\colon M\to K$ is a $k$-sheeted covering up to the boundary, and every point of $C_0$ has exactly $k$ preimages in $\partial M$. Since $x|_{\partial M}$ is injective, $k=1$, and $x$ is a diffeomorphism onto $K$. Finally, $\rho_i=\operatorname{Per}(C_i)/(2\pi)$ is the radius of $C_i$.
\end{proof}

For $\rho_0=\rho_1$ the extremal piece is the marginally stable catenoid piece between two equal coaxial circles. Theorem~\ref{thm:rigidity} is the analogue of the rigidity statement in \cite[Thm.~4.1, Prop.~4.4]{BernsteinBreiner14}. It applies to boundaries with several components and to non-convex curves, and it forces such boundaries to be two coaxial circles. The injectivity of $x|_{\partial M}$ cannot be dropped: a $k$-fold covering of the extremal catenoid piece also attains equality. It can be weakened to the assumption that some boundary component is mapped to $C_0$ with degree $\pm1$. The regularity of $x$ up to the boundary is used only in the last paragraph of the proof, and we expect that it can be removed by a fundamental-group argument.

It is useful to compare our result with two related nonexistence results in the literature, which are formulated in terms of different geometric data.

\begin{enumerate}[label=(\arabic*)]
\item \emph{Bernstein--Breiner.} Using the Osserman--Schiffer inequality \[L''(t)\ge L(t)\] for the lengths of the images of the circles $|z|=e^t$ under a conformal harmonic annulus \cite{OssermanSchiffer75}, \cite[Lemma~2.3]{BernsteinBreiner14}, Bernstein and Breiner proved a sharp nonexistence result in terms of the \emph{lengths} $|\Gamma_i|$ \cite[Thm.~4.1, Prop.~4.4]{BernsteinBreiner14}. It holds for connected immersed surfaces of arbitrary genus, and its equality case is rigid. Their threshold $F_\Omega$ is characterized in \cite[Lemma~4.2(1)]{BernsteinBreiner14}: there is a vertical catenoid piece in the slab $\Omega$ with boundary lengths $L_-,L_+$ if and only if $L_+\ge F_\Omega(L_-)$, and \cite[Prop.~4.4]{BernsteinBreiner14} excludes $|\Gamma_+|<F_\Omega(|\Gamma_-|)$.

The two thresholds coincide. The heights of vertical catenoid pieces bounded by coaxial circles of radii $\rho_0,\rho_1$ fill the interval $(0,D(\rho_0,\rho_1)]$: as the waist $m$ runs over $(0,\min\{\rho_0,\rho_1\}]$, the height \[m[\operatorname{arccosh}(\rho_0/m)+\operatorname{arccosh}(\rho_1/m)]\] is continuous and tends to $0$ as $m\to0$, and pieces not containing the neck are lower. Hence $\tau_1-\tau_0\le D(\rho_0,\rho_1)$ if and only if such a catenoid piece fits into the slab, i.e. if and only if $2\pi\rho_1\ge F_\Omega(2\pi\rho_0)$. 
Therefore Corollary~\ref{cor:BB} is the nonexistence part of \cite[Prop.~4.4]{BernsteinBreiner14} with $|\Gamma_i|$ replaced by $\operatorname{Per}(\conv\Gamma_i)$. Since $\operatorname{Per}(\conv\Gamma)\le|\Gamma|$, with equality exactly for convex $\Gamma$, and $D$ is increasing in each argument:
\begin{itemize}
\item for convex curves, Corollary~\ref{cor:BB} and Theorem~\ref{thm:rigidity} reprove their result, including rigidity, by a different method;
\item for non-convex curves, Corollary~\ref{cor:BB} is strictly stronger than their nonexistence result, since $\operatorname{Per}(\conv\Gamma_i)<|\Gamma_i|$.
\end{itemize}
\item \emph{Diameter bounds.} For compact minimal surfaces bounded by one curve in each of two parallel planes, Rossman \cite{Rossman08} proved that the distance $h$ between the planes satisfies $h\le\max\{d_0,d_1\}$, where $d_i$ is the diameter of $\Gamma_i$. This generalizes an earlier bound of Nitsche for annuli with the constant $\frac32$, see \cite{Rossman08}. Corollary~\ref{cor:BB} improves the constant. By Cauchy's formula 
    \[\operatorname{Per}(K)=\int_0^\pi w_K(\varphi)\,d\varphi,\] where $w_K\le\operatorname{diam}K$ is the width, we have $\operatorname{Per}(\conv\Gamma_i)\le\pi d_i$, i.e. $\rho_i\le d_i/2$. Since $D$ is increasing in each argument, we obtain (see \eqref{xnoll})
\begin{align*}
h\ &\le\ D(\rho_0,\rho_1)\ \le\ D(\rho_{\max},\rho_{\max})=c_*\rho_{\max}\\
&\le\ \tfrac{c_*}{2}\max\{d_0,d_1\}=\frac{\max\{d_0,d_1\}}{\sinh x_0}\approx0.6627\,\max\{d_0,d_1\}.
\end{align*}
\item \emph{Intermediate sections.} Corollary~\ref{cor:BB} is a consequence of the bound 
\[B(\tau)\ge B_{\min}\cosh((\tau-\tau_m)/B_{\min}),\] 
which controls the hull of \emph{every} intermediate section, not only the boundary.
\item \emph{Lifetime of tubes.} For $n\ge3$ and $p=2$, Vedenyapin--Miklyukov \cite{MV86} and Miklyukov--Tkachev \cite{MT89} showed that minimal tubes have bounded lifetime, using the distance to an axis; see also Klyachin \cite{Klyachin92} for arbitrary codimension. In \cite[Cor.~5]{T97} this was extended to $p$-minimal tubes with $n>p$, using the circumradius $r(M)=\min R$, under the condition $\sigma\ge\varepsilon$. Theorem~\ref{thm:integrated}(i) replaces $r(M)$ by the smaller quantity $B_{\min}$ and removes the condition $\sigma\ge\varepsilon$: by Theorem~\ref{thm:semiconvex}, the joint regularity which this condition was meant to provide is automatic.
\item \emph{Rigidity.} Theorem~\ref{thm:rigidity} shows that equality in Corollary~\ref{cor:BB} holds only for the extremal catenoid piece. The proof uses the comparison principle only through the integrated inequality: equality forces the hulls to be coaxial discs. The extreme-point property \ref{S5} then places the whole catenoid inside $M$, and analyticity does the rest. No equality case of Theorem~\ref{thm:main} at individual contact points is needed.
\end{enumerate}

\begin{remark}[Edges shorten the surface]\label{rem:strict}
Suppose that $\lambda>0$ on a subset of $(\tau_0,\tau_1)$ of positive measure; for $n=2$ this holds if the hulls $\Omega(\tau)$ have edges for a set of levels of positive measure. Then the inequalities in Theorem~\ref{thm:integrated}(i), (ii) and in Corollary~\ref{cor:BB} are strict. Indeed, in the proof of Lemma~\ref{lem:ODE}, \eqref{eq:Bineqlambda} improves the estimate of $d\mu$ to
\[
d\mu\ \le\ -\,B'\,\Phi(B')\,\frac{\lambda}{B-\lambda}\,d\tau\qquad\text{on }(\tau_m,\tau_1).
\]
Assume that $\{\lambda>0\}\cap(\tau_m,\tau_1)$ has positive measure. Since $B'>0$ on $(\tau_m,\tau_1)$, there is $\tau'<\tau_1$ with $\mu<B_{\min}$ on $(\tau',\tau_1)$. As $\mu=B\,\Phi(B')$ means $B'=\sqrt{(B/\mu)^{2\beta}-1}$, we get $B'>\sqrt{(B/B_{\min})^{2\beta}-1}$ on $(\tau',\tau_1)$, and the bound $\tau_1-\tau_m\le B_{\min}G_\beta(B_1/B_{\min})$ becomes strict. The interval $(\tau_0,\tau_m)$ is treated in the same way.
\end{remark}

\section{Application to the Meeks conjecture}\label{sec:MAC}

In this section $p=2$, $n=2$, and $M\subset\R^3$ is a compact connected embedded minimal surface with $\partial M=\Gamma_0\cup\Gamma_1$, $\Gamma_i\subset\Pi_{\tau_i}$ Jordan curves, $\tau_0<\tau_1$. If $\Gamma_0,\Gamma_1$ are convex, embeddedness is automatic by \cite{EWW02}. Since $x_3$ is harmonic and nonconstant, the strong maximum principle \cite[Thm.~3.5]{GT01} gives $x_3(M^\circ)\subset(\tau_0,\tau_1)$.

\begin{lemma}\label{lem:extreme}
For $\tau\in(\tau_0,\tau_1)$, the set $\partial\Omega(\tau)\setminus\Sigma(\tau)$ is a union of open segments, each contained in an edge of $\Omega(\tau)$ and with both endpoints in $\Sigma(\tau)$.
\end{lemma}

\begin{proof}
By \ref{S5}, every extreme point of $\Omega(\tau)$ lies in $\Sigma(\tau)$. A point $z\in\partial\Omega(\tau)\setminus\Sigma(\tau)$ is therefore not extreme, so it lies in the relative interior of an edge $E\subset\partial\Omega(\tau)$, whose endpoints are extreme. The component of $E\setminus\Sigma(\tau)$ containing $z$ is an open segment, and its endpoints lie in the closed set $\Sigma(\tau)$.
\end{proof}

\begin{definition}\label{def:bridge}
A \emph{supporting chord} at level $\tau$, briefly a \emph{chord}, is a segment $[a,b]\subset\partial\Omega(\tau)$, $a\ne b$, with $a,b\in\Sigma(\tau)$ and $(a,b)\cap\Sigma(\tau)=\emptyset$. The open segment $(a,b)$ is called a \emph{gap}.\footnote{In computational geometry, the common supporting edge of the convex hulls of two point sets is called a \emph{bridge} \cite{KirkpatrickSeidel86}; this is the situation of Figure~\ref{fig:bridge}(b). We do not use that term, to avoid confusion with the bridge principle in minimal surface theory.}
\end{definition}

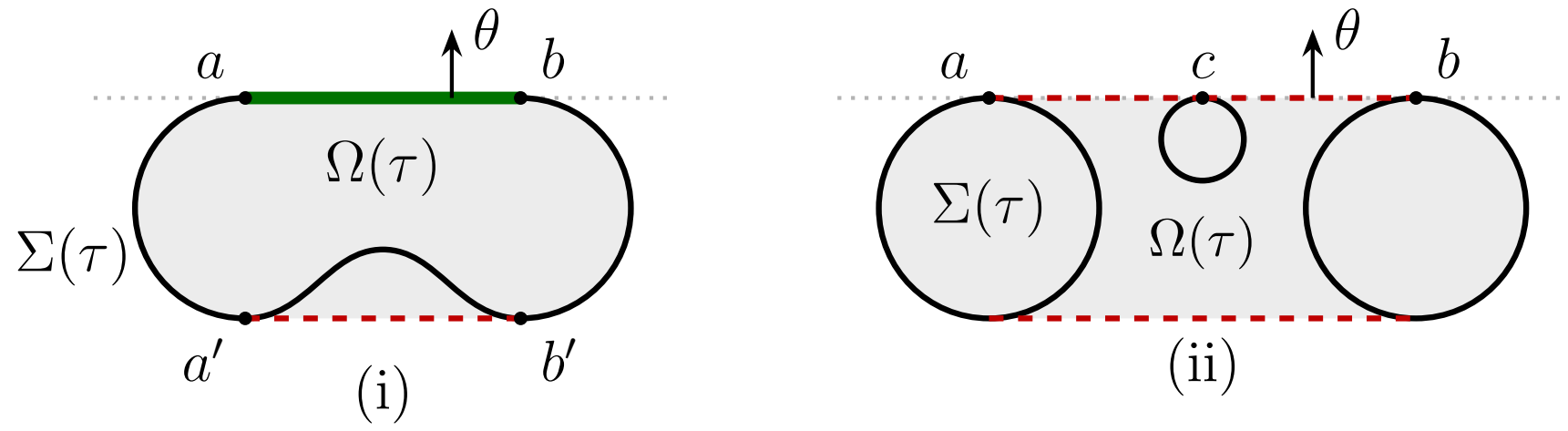
\begin{figure}[ht]
\centering
\begin{tikzpicture}[scale=1.0,
  sect/.style={very thick, black},
  hull/.style={fill=gray!15},
  chord/.style={very thick, dashed, red!75!black},
  inside/.style={line width=2.6pt, green!45!black}]
\begin{scope}
 \fill[hull] (-1,0.8) -- (1,0.8) arc[start angle=90,end angle=-90,radius=0.8]
   -- (-1,-0.8) arc[start angle=270,end angle=90,radius=0.8] -- cycle;
 \draw[gray!60, dotted, thick] (-2.1,0.8) -- (2.1,0.8);
 \draw[sect] (-1,0.8) -- (1,0.8) arc[start angle=90,end angle=-90,radius=0.8]
   .. controls (0.55,-0.8) and (0.4,-0.3) .. (0,-0.3)
   .. controls (-0.4,-0.3) and (-0.55,-0.8) .. (-1,-0.8)
   arc[start angle=270,end angle=90,radius=0.8];
 \draw[inside] (-1,0.8) -- (1,0.8);
 \draw[chord] (-1,-0.8) -- (1,-0.8);
 \fill (-1,0.8) circle (1.4pt) node[above left] {$a$};
 \fill (1,0.8) circle (1.4pt) node[above right] {$b$};
 \fill (-1,-0.8) circle (1.4pt) node[below left] {$a'$};
 \fill (1,-0.8) circle (1.4pt) node[below right] {$b'$};
 \draw[-{Stealth[length=2mm]}, thick] (0.5,0.8) -- (0.5,1.3) node[right] {$\theta$};
 \node at (0,0.3) {$\Omega(\tau)$};
 \node at (-2.25,-0.35) {$\Sigma(\tau)$};
 \node at (0,-1.35) {(i)};
\end{scope}
\begin{scope}[xshift=5.6cm]
 \fill[hull] (-1.2,0.8) -- (1.9,0.8) arc[start angle=90,end angle=-90,radius=0.8]
      -- (-1.2,-0.8) arc[start angle=270,end angle=90,radius=0.8] -- cycle;
 \draw[gray!60, dotted, thick] (-2.3,0.8) -- (3.0,0.8);
 \draw[sect] (-1.2,0) circle (0.8);
 \draw[sect] (1.9,0) circle (0.8);
 \draw[sect] (0.35,0.5) circle (0.3);
 \draw[chord] (-1.2,0.8) -- (0.35,0.8);
 \draw[chord] (0.35,0.8) -- (1.9,0.8);
 \draw[chord] (-1.2,-0.8) -- (1.9,-0.8);
 \fill (-1.2,0.8) circle (1.4pt) node[above left] {$a$};
 \fill (0.35,0.8) circle (1.4pt) node[above] {$c$};
 \fill (1.9,0.8) circle (1.4pt) node[above right] {$b$};
 \draw[-{Stealth[length=2mm]}, thick] (1.15,0.8) -- (1.15,1.3) node[right] {$\theta$};
 \node at (-1.2,0) {$\Sigma(\tau)$};
 \node[font=\small] at (0.35,-0.25) {$\Omega(\tau)$};
 \node at (0.35,-1.15) {(ii)};
\end{scope}
\end{tikzpicture}
\caption{Edges of $\Omega(\tau)$ which are not chords; chords are dashed. (i) The edge $[a,b]$ lies in $\Sigma(\tau)$ (thick); the chord $[a',b']$ across the dent is a chord. By Lemma~\ref{lem:nosegment}, configuration (i) is impossible for minimal surfaces. (ii) The edge $[a,b]$ meets $\Sigma(\tau)$ at an interior point $c$; it is the union of the two chords $[a,c]$ and $[c,b]$. Such triple contacts are non-generic (Remark~\ref{rem:generic}). In both cases all points of $\partial\Omega(\tau)\cap\Sigma(\tau)$, including $c$, are regular points of $x_3$ by Lemma~\ref{lem:regular}, and $\Omega(\tau)$ has no corners, as required by Lemma~\ref{lem:smooth}.}
\label{fig:nonbridge}
\end{figure}

By Lemma~\ref{lem:extreme}, the gaps are precisely the components of $\partial\Omega(\tau)\setminus\Sigma(\tau)$, so condition \eqref{eq:F} fails if and only if some intermediate section has a chord. An edge of $\Omega(\tau)$ need not itself be a chord: it may lie in $\Sigma(\tau)$, or it may meet $\Sigma(\tau)$ at interior points, in which case it is a union of several chords (Figure~\ref{fig:nonbridge}). The first possibility does not occur for minimal surfaces.

\begin{lemma}\label{lem:nosegment}
For every $\tau\in(\tau_0,\tau_1)$, the section $\Sigma(\tau)$ contains no nondegenerate segment, and its intersection with any straight line is finite.
\end{lemma}

\begin{proof}
In conformal coordinates the immersion $x$ is harmonic, hence real-analytic \cite{Osserman86,DHS10}. Let $m\in M^\circ$ and let $z$ be a conformal coordinate centred at $m$. If $v$ is a nonconstant harmonic function near $m$ with $v(m)=0$, then $v=\mathrm{Re}\,F$ with $F$ holomorphic, $F(0)=0$. If $F$ vanishes to order $k+1\ge1$ at $0$, then $F=\phi^{k+1}$ with $\phi$ biholomorphic near $0$. Hence $\{v=0\}$ near $m$ is the union of $k+1$ real-analytic embedded arcs through $m$, the preimages under $\phi$ of the lines $\mathrm{Re}\,w^{k+1}=0$.

Suppose that $\Sigma(\tau)$ contains a segment $[a,b]$ with $a\ne b$. Put $u=(b-a)/|b-a|$ and $\gamma(t)=a+tu$, and let
\[
T=\sup\{t\ge0\mid \gamma([0,t])\subset\Sigma(\tau)\}
\]
be the length of the longest segment in $\Sigma(\tau)$ starting at $a$ in the direction $u$. Then $|b-a|\le T<\infty$ by compactness, and $\gamma([0,T])\subset\Sigma(\tau)$ since $\Sigma(\tau)$ is closed. We show, by analytic continuation, that the segment can be prolonged beyond $\gamma(T)$; this contradicts the choice of $T$.

\emph{The endpoint.} Let $m$ be the point with $x(m)=(\gamma(T),\tau)$. It is an interior point of $M$, since $\tau\in(\tau_0,\tau_1)$. Apply the local description above to the harmonic function $v=x_3-\tau$: near $m$, the level set $\{v=0\}$ is a finite union of real-analytic arcs through $m$.

\emph{The segment lies on one analytic arc.} Since $x$ is an embedding, for small $\delta>0$ the preimage $x^{-1}\big(\gamma((T-\delta,T))\times\{\tau\}\big)$ is a connected arc in $\{v=0\}\setminus\{m\}$ ending at $m$. Hence it is contained in one of these analytic arcs, say $\beta\colon(-\varepsilon,\varepsilon)\to M^\circ$ with $\beta(0)=m$, and corresponds to the parameters $s\in(-\varepsilon,0)$.

\emph{Analytic continuation and contradiction.} Choose a unit vector $\vartheta\perp u$. The real-analytic functions
\[
s\mapsto x_3(\beta(s))-\tau,\qquad s\mapsto\scal{y(\beta(s))-a}{\vartheta}
\]
vanish on $(-\varepsilon,0)$, hence on all of $(-\varepsilon,\varepsilon)$. So $x(\beta(s))$ lies on the line $\gamma(\R)\times\{\tau\}$ for every $s\in(-\varepsilon,\varepsilon)$.
Since $x$ is injective, the function $s\mapsto\scal{y(\beta(s))-a}{u}$ is continuous and injective, hence strictly monotone. It is less than $T$ for $s<0$ and equal to $T$ at $s=0$, so it is greater than $T$ for $s>0$. Thus $\gamma([T,T+\varepsilon'))\subset\Sigma(\tau)$ for some $\varepsilon'>0$, contradicting the definition of $T$.

For the second statement, $x_3^{-1}(\tau)$ is compact and, by the local description above, a finite union of compact real-analytic arcs. On each arc a linear function $\scal{y}{\vartheta}-c$ either vanishes identically, which would produce a segment in $\Sigma(\tau)$, or has finitely many zeros.
\end{proof}

\begin{corollary}\label{cor:C1}
For every $\tau\in(\tau_0,\tau_1)$ each edge of $\Omega(\tau)$ is the union of finitely many chords. Consequently, the following are equivalent:
\begin{enumerate}[label=(\alph*)]
\item condition \eqref{eq:F} holds;
\item no intermediate section has a chord;
\item $H(\tau,\cdot)\in C^1(\R^2\setminus0)$ for every $\tau\in(\tau_0,\tau_1)$;
\item the area function $S_\tau=S(\Omega(\tau),\cdot)$ has no atoms for every $\tau\in(\tau_0,\tau_1)$.
\end{enumerate}
Moreover, the atomic part of $S_\tau$ is $\sum_b|b|\,\delta_{\theta_b}$, where $b$ runs over the chords at level $\tau$ and $\theta_b$ is the outer normal of $\Omega(\tau)$ along $b$.
\end{corollary}

\begin{proof}
Let $E$ be an edge of $\Omega(\tau)$. Its endpoints are extreme, hence lie in $\Sigma(\tau)$ by \ref{S5}. By Lemma~\ref{lem:nosegment}, $E\cap\Sigma(\tau)$ is a finite set containing the endpoints, so $E$ is the union of the finitely many chords between consecutive points of $E\cap\Sigma(\tau)$. The equivalence (a)$\Leftrightarrow$(b) was noted above. By \ref{S3}, (c) holds if and only if every face of every $\Omega(\tau)$ is a point, i.e. $\Omega(\tau)$ has no edges, which is equivalent to (b) by the first assertion. By \ref{S6}, the atoms of $S_\tau$ are the edges $E$ with masses $|E|$, which gives (c)$\Leftrightarrow$(d) and, since each edge is the union of its chords, the last assertion.
\end{proof}

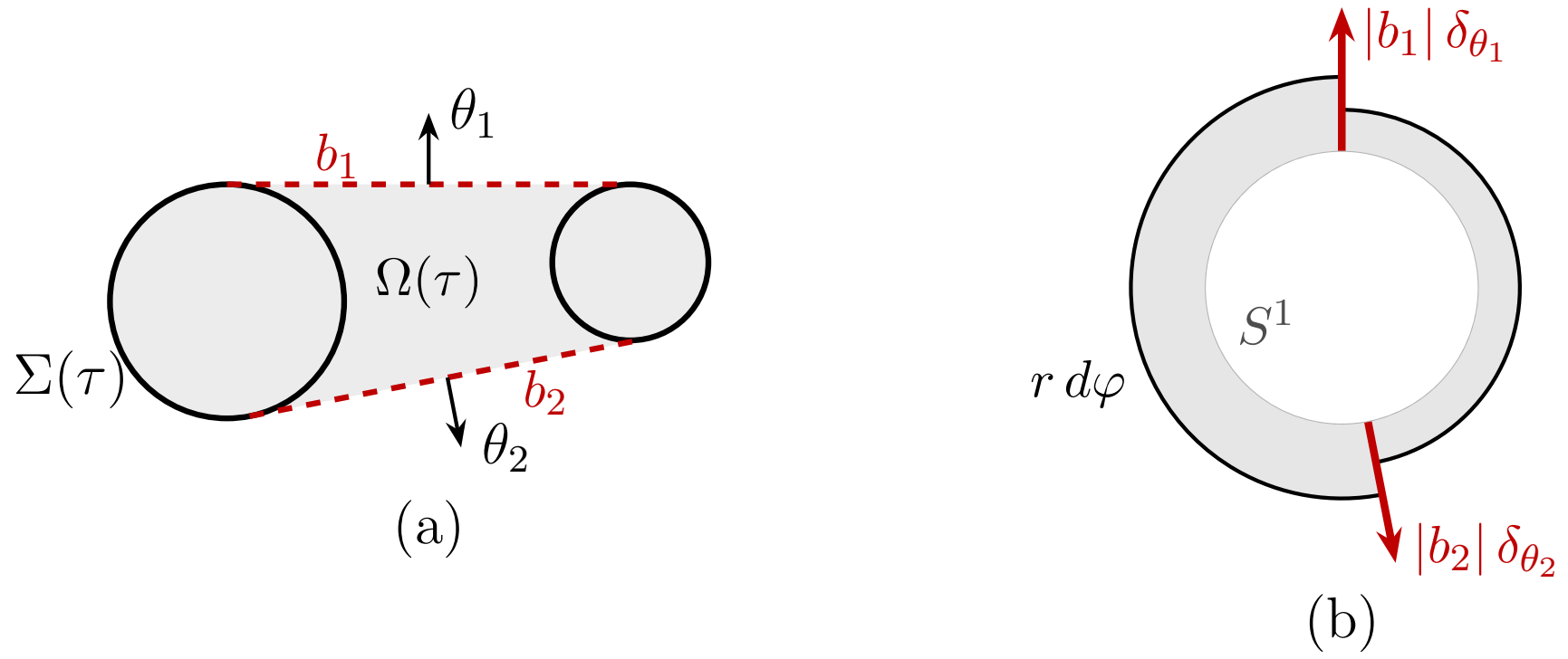
\begin{figure}[ht]
\centering
\begin{tikzpicture}[scale=1.0,
  sect/.style={very thick, black},
  hull/.style={fill=gray!15},
  chord/.style={very thick, dashed, red!75!black}]
\begin{scope}[scale=0.95]
 \fill[hull] (-1.2,0.9) arc[start angle=90,end angle=281.06,radius=0.9]
     -- (2.015,-0.289) arc[start angle=-78.94,end angle=90,radius=0.6] -- cycle;
 \draw[sect] (-1.2,0) circle (0.9);
 \draw[sect] (1.9,0.3) circle (0.6);
 \draw[chord] (-1.2,0.9) -- (1.9,0.9);
 \draw[chord] (-1.027,-0.883) -- (2.015,-0.289);
 \draw[-{Stealth[length=2mm]}, thick] (0.35,0.9) -- (0.35,1.45) node[right] {$\theta_1$};
 \draw[-{Stealth[length=2mm]}, thick] (0.494,-0.586) -- ++(0.105,-0.539) node[right] {$\theta_2$};
 \node[red!75!black, font=\small] at (-0.35,1.12) {$b_1$};
 \node[red!75!black, font=\small] at (1.25,-0.7) {$b_2$};
 \node[font=\small] at (0.35,0.15) {$\Omega(\tau)$};
 \node at (-2.4,-0.6) {$\Sigma(\tau)$};
 \node at (0.35,-1.75) {(a)};
\end{scope}
\begin{scope}[xshift=7.0cm, yshift=0.1cm]
 \draw[gray!60] (0,0) circle (1.0);
 \node[gray!60!black, font=\small] at (-0.55,-0.25) {$S^1$};
 \fill[gray!20] (90:1.54) arc[start angle=90,end angle=281.06,radius=1.54]
   -- (281.06:1.0) arc[start angle=281.06,end angle=90,radius=1.0] -- cycle;
 \fill[gray!20] (-78.94:1.30) arc[start angle=-78.94,end angle=90,radius=1.30]
   -- (90:1.0) arc[start angle=90,end angle=-78.94,radius=1.0] -- cycle;
 \draw[thick] (90:1.54) arc[start angle=90,end angle=281.06,radius=1.54];
 \draw[thick] (-78.94:1.30) arc[start angle=-78.94,end angle=90,radius=1.30];
 \draw[red!75!black, line width=1.6pt, -{Stealth[length=2.5mm]}] (90:1.0) -- (90:2.05);
 \draw[red!75!black, line width=1.6pt, -{Stealth[length=2.5mm]}] (-78.94:1.0) -- (-78.94:2.05);
 \node[red!75!black, font=\small, right] at (90:1.85) {$|b_1|\,\delta_{\theta_1}$};
 \node[red!75!black, font=\small, right] at (-78.94:1.95) {$|b_2|\,\delta_{\theta_2}$};
 \node[font=\small] at (200:2.05) {$r\,d\varphi$};
 \node at (0,-2.45) {(b)};
\end{scope}
\end{tikzpicture}
\caption{(a) A section with two components; its hull has the two chords $b_1,b_2$ with outer normals $\theta_1,\theta_2$. (b) The area function $S_\tau=S(\Omega(\tau),\cdot)$ on $S^1$: an absolutely continuous part with density $r=h_{\varphi\varphi}+h$, the radius of curvature of the arc of $\Sigma(\tau)$ with normal $\theta(\varphi)$ (shaded, drawn radially), and one atom $|b|\,\delta_{\theta_b}$ per chord. $S_\tau$ charges every open arc because $\Omega(\tau)$ has no corners (Lemma~\ref{lem:smooth}). By Corollary~\ref{cor:C1}, MAC is equivalent to the absence of atoms at all levels.}
\label{fig:measure}
\end{figure}

\begin{remark}\label{rem:generic}
(1) Configuration (ii) of Figure~\ref{fig:nonbridge} is not excluded by analyticity, but it is non-generic. For fixed $\tau$, supporting lines touching $\Sigma(\tau)$ at two points occur only in isolated directions $\theta$, and a third contact point on the same line is one more equation; so in a generic one-parameter family of sections it occurs only at isolated levels $\tau$.

(2) For $p\ne2$, the hypersurface is real-analytic on $\Mreg$ by Morrey's theorem \cite{Morrey58}, since there it satisfies an elliptic equation with analytic coefficients, cf.\ \cite[(eq)]{T97}. The argument of Lemma~\ref{lem:nosegment} therefore applies as long as the continued segment stays in $\Mreg$. It may fail only if the segment runs into a critical point of $x_{n+1}$, where the $p$-Laplacian degenerates; by \cite[Lemma~1]{T97} such points are flat points of the hypersurface.
\end{remark}

\begin{theorem}\label{thm:MAC}
Let $p=2$, $n=2$, and let $M\subset\R^3$ be a compact connected embedded minimal surface with $\partial M=\Gamma_0\cup\Gamma_1$, $\Gamma_i\subset\Pi_{\tau_i}$ Jordan curves, $\tau_0<\tau_1$. Then:
\begin{enumerate}[label=(\roman*)]
\item If \eqref{eq:F} holds, then $x_3$ has no critical points on $M^\circ$, every section $\Sigma(\tau)$, $\tau\in(\tau_0,\tau_1)$, is the convex Jordan curve $\partial\Omega(\tau)$, and $M$ is an annulus.
\item If $\Gamma_0,\Gamma_1$ are convex and $M$ is an annulus, then \eqref{eq:F} holds.
\end{enumerate}
Consequently, for convex $\Gamma_0,\Gamma_1$, the Meeks annulus conjecture holds for $M$ if and only if no intermediate section of $M$ has a chord, if and only if all support functions $H(\tau,\cdot)$, $\tau\in(\tau_0,\tau_1)$, are of class $C^1$.
\end{theorem}

\begin{proof}
(ii) By Shiffman's theorem \cite{Shiffman56}, every $\Sigma(\tau)$ is a convex Jordan curve, hence equals $\partial\Omega(\tau)$.

(i) Fix $\tau\in(\tau_0,\tau_1)$ and let $\mathcal C(\tau)=x^{-1}(\partial\Omega(\tau)\times\{\tau\})\subset M^\circ$. By \eqref{eq:F} and the definition of contact points, every point of $\mathcal C(\tau)$ is a contact point, hence regular by Lemma~\ref{lem:regular}. Near a regular point $m$, the level set $x_3^{-1}(\tau)$ is an embedded $C^2$ arc. Since $x$ is an embedding of the compact $M$, small balls around $x(m)$ meet $x(M)$ only in the image of a small neighbourhood of $m$.

\emph{Step 1: $\Omega(\tau)$ has interior points.} This is Lemma~\ref{lem:smooth}.

\emph{Step 2: $\mathcal C(\tau)$ is a component of $x_3^{-1}(\tau)$.} By Step 1, $\partial\Omega(\tau)$ is a Jordan curve. Near each of its points it is contained in the regular arc $\Sigma(\tau)$ through that point, so by invariance of domain in dimension one it contains a neighbourhood of that point in $\Sigma(\tau)$. Thus $\mathcal C(\tau)$ is open and closed in $x_3^{-1}(\tau)$, and it is connected.

\emph{Step 3: $\mathcal E=\bigcup_{\tau}\mathcal C(\tau)$ is open and closed in $M^\circ$.}

We first prove that $\mathcal E$ is closed in $M^\circ$. Let $m_k\in\mathcal E$ converge to $m\in M^\circ$.
Then $\scal{y(m_k)}{\theta_k}=H(x_3(m_k),\theta_k)$ for some $\theta_k\in S^1$. Passing to a subsequence $\theta_k\to\theta$ and using the continuity of $H$ (Lemma~\ref{lem:convex}), we get $\scal{y(m)}{\theta}=H(x_3(m),\theta)$. Hence $y(m)\in\partial\Omega(x_3(m))$, that is, $m\in\mathcal E$.

We now turn to openness. Fix $\tau$. We show that $\mathcal C(\tau')\subset\mathcal E$ fills a whole neighbourhood of $\mathcal C(\tau)$ for $\tau'$ near $\tau$. This needs a neighbourhood of $\mathcal C(\tau)$ in which the level sets of $x_3$ are explicitly known. Since $\mathcal C(\tau)$ is a compact curve of regular points, the standard collar construction gives one (cf.\ the proof of Theorem~3.1 in \cite[\S3]{Milnor63}). The function $|\nabla x_3|$ is bounded below by a positive constant near $\mathcal C(\tau)$, so the vector field $X=\nabla x_3/|\nabla x_3|^2$ is defined there. It satisfies $X(x_3)=1$, so its flow $\Phi_t$ maps the level $x_3^{-1}(\tau)$ to the level $x_3^{-1}(\tau+t)$. For small $\varepsilon>0$, the map $(q,t)\mapsto\Phi_t(q)$ is therefore a diffeomorphism of $\mathcal C(\tau)\times(-\varepsilon,\varepsilon)$ onto an open neighbourhood $V$ of $\mathcal C(\tau)$. In $V$ every level set is a single closed curve:
\[
V\cap x_3^{-1}(\tau')=\Phi_{\tau'-\tau}\big(\mathcal C(\tau)\big),\qquad |\tau'-\tau|<\varepsilon.
\]

It remains to show that $\mathcal C(\tau')$ is exactly this curve. Since $x$ is an embedding, $x(M\setminus V)$ is compact and disjoint from $x(\mathcal C(\tau))$, hence at positive distance $d$ from it. By continuity of $H$ and Step 1, $\partial\Omega(\tau')\to\partial\Omega(\tau)$ in the Hausdorff distance \cite[\S1.8]{Schneider14}. So for $\tau'$ close to $\tau$, every point of $x(\mathcal C(\tau'))$ lies within distance $d$ of $x(\mathcal C(\tau))$, and hence $\mathcal C(\tau')\subset V$. By Step 2 at level $\tau'$, $\mathcal C(\tau')$ is a component of $x_3^{-1}(\tau')$; it is nonempty, since the extreme points of $\Omega(\tau')$ lie in $\Sigma(\tau')$ by \ref{S5}. Since it is contained in the connected curve $V\cap x_3^{-1}(\tau')$, the two coincide. Thus $V\subset\mathcal E$.

\emph{Step 4: conclusion.} Since $M^\circ$ is connected and $\mathcal E\neq\emptyset$, we have
$\mathcal E=M^\circ$. Hence every level set $x_3^{-1}(\tau)$ is the single
regular curve $\mathcal C(\tau)$, and $x_3$ has no critical points on
$M^\circ$. The map
\[
x_3:M^\circ\to(\tau_0,\tau_1)
\]
is therefore a proper submersion with fibre $S^1$. By Ehresmann's theorem
\cite{Ehresmann51} it is a locally trivial fibration, and, since
$(\tau_0,\tau_1)$ is contractible, the fibration is trivial. Thus
\[
M^\circ\cong S^1\times(\tau_0,\tau_1).
\]
Since $\partial M=\Gamma_0\sqcup\Gamma_1$ consists of two Jordan curves,
it follows that $M$ is an annulus.

\end{proof}

The following proposition shows that the endpoints of a chord are highly constrained.

\begin{proposition}[Chords]\label{prop:bridge}
Let $[a,b]$ be a chord at level $\tau$, contained in the edge of $\Omega(\tau)$ with outer normal $\theta$, and let $m_a,m_b\in M^\circ$ be the points with $x(m_a)=(a,\tau)$, $x(m_b)=(b,\tau)$. Then:
\begin{enumerate}[label=(\roman*)]
\item $m_a,m_b$ are regular points of $x_3$;
\item $x(M)$ is disjoint from the open convex set $K_\theta=\{(y,t):\ \scal{y}{\theta}>H(t,\theta)\}$;
\item if $t\mapsto H(t,\theta)$ is differentiable at $\tau$, which holds for all but countably many $\tau$, then
\[
T_{m_a}M=T_{m_b}M=\Pi:=\{(y,t):\ \scal{y-a}{\theta}=H_\tau(\tau,\theta)\,(t-\tau)\},
\]
and $\Pi$ supports $\overline{K_\theta}$ along the line through $a$ and $b$.
\end{enumerate}
\end{proposition}

\begin{proof}
(i) is Lemma~\ref{lem:regular}. (ii) is \eqref{eq:basic}; $K_\theta$ is convex as a strict superlevel set of the concave function $\scal{y}{\theta}-H(t,\theta)$ (Lemma~\ref{lem:convex}). For (iii), the function $\psi=\scal{y}{\theta}-H(x_3,\theta)$ attains its maximum $0$ at $m_a$ and is differentiable there, so $\nabla\psi(m_a)=\big((\theta,0)-H_\tau e\big)^\top=0$. Thus $(\theta,-H_\tau)$ is normal to $M$ at $m_a$, and likewise at $m_b$. The plane through $(a,\tau)$ with this normal contains $(b,\tau)$, since $\scal{b-a}{\theta}=0$. By convexity, $H(t,\theta)\ge H(\tau,\theta)+H_\tau(t-\tau)$, so $K_\theta$ lies on one side of $\Pi$, and $\Pi\cap\overline{K_\theta}$ contains the line through $a,b$. A convex function of one variable is differentiable outside a countable set.
\end{proof}

\begin{remark}
At the remaining countably many levels, Lemma~\ref{lem:smooth} still gives $\nu(m_a)\parallel(\theta,-c_a)$ and $\nu(m_b)\parallel(\theta,-c_b)$ with $c_a,c_b\in[\partial_\tau^-H(\tau,\theta),\partial_\tau^+H(\tau,\theta)]$. So the tangent planes at $m_a$ and $m_b$ both contain the line through $a$ and $b$, but they may differ.
\end{remark}

\begin{remark}[Towards excluding chords]\label{rem:program}
Proposition~\ref{prop:bridge} shows that a \textit{hypothetical counterexample} to the Meeks annulus conjecture contains a configuration in which $M$ is tangent to one plane $\Pi$ at two points joined by a horizontal segment whose interior misses $M$, while $M$ avoids the convex set $K_\theta$ supported by $\Pi$ along that segment. When $\Gamma_0,\Gamma_1$ are strictly convex, $\Pi$ meets $\partial M$ in at most four points. Combined with the structure of $M\cap\Pi$ near the two tangency points, which consists of at least two analytic arcs crossing at each, i.e.\ at least four half-branches (apply the local description in the proof of Lemma~\ref{lem:nosegment} to $v=\scal{x-x(m_a)}{\nu(m_a)}$), this suggests a Rad\'o-type counting argument. A second route is a comparison principle for \eqref{eq:MA} against the stable annulus of \cite{MeeksWhite91}. Both are left for further work, beginning with the case of two circles.
\end{remark}

\section{Several boundary curves}\label{sec:several}

We now allow each boundary set to have several connected components. The comparison principle of Sections~\ref{sec:contact}--\ref{sec:integrated} does not require $\Gamma_0$ or $\Gamma_1$ to be connected or convex, and therefore extends directly to this setting.
In this section $M\subset\R^3$ is a compact connected immersed minimal surface, smooth up to the boundary, with
\[
\partial M=\Gamma_0\cup\Gamma_1,\qquad \Gamma_i\subset\Pi_{\tau_i},
\]
where $\Gamma_i$ consists of $k_i\ge1$ pairwise disjoint smooth closed curves; we put $d=\tau_1-\tau_0$. As in Section~\ref{sec:MAC}, $x_3(M^\circ)\subset(\tau_0,\tau_1)$ by the strong maximum principle, and $\nabla x_3\neq0$ on $\partial M$ by the Hopf boundary lemma \cite[Lemma~3.4]{GT01}.

Theorem~\ref{thm:integrated} and its proof apply without change.

\begin{corollary}\label{cor:BBseveral}
Let $\rho_i=\operatorname{Per}(\conv\Gamma_i)/(2\pi)$. Then $d\le D(\rho_0,\rho_1)$.
\end{corollary}

The nonexistence theorem of Bernstein and Breiner \cite[Thm.~4.1, Prop.~4.4]{BernsteinBreiner14} assumes that each of the two boundary sets is a single closed curve. For several curves, the total length $|\Gamma_i|$ would in any case be the wrong quantity: it does not see how the curves are placed relative to each other. The hull perimeter does. We now show that the inequality is strict as soon as no single boundary curve carries the hull of its plane.

\begin{lemma}\label{lem:edgesnear}
Suppose that $\conv\gamma\ne\conv\Gamma_0$ for every component $\gamma$ of $\Gamma_0$. Then there is $\varepsilon>0$ such that $\Omega(\tau)$ has an edge for every $\tau\in(\tau_0,\tau_0+\varepsilon)$.
\end{lemma}

\begin{proof}
Since $\nabla x_3\ne0$ on $\Gamma_0$, there is $\varepsilon>0$ such that the flow of $\nabla x_3/|\nabla x_3|^2$ identifies $x_3^{-1}([\tau_0,\tau_0+\varepsilon))$ with $\Gamma_0\times[\tau_0,\tau_0+\varepsilon)$. For $\tau$ in this range,
\[
\Sigma(\tau)=C_1(\tau)\cup\dots\cup C_{k_0}(\tau),
\]
where $C_j(\tau)$ is the projection of the flow image of the $j$-th component $\gamma_j$ of $\Gamma_0$, and $C_j(\tau)\to\gamma_j$ in the Hausdorff distance as $\tau\downarrow\tau_0$. Since the $\gamma_j$ are pairwise disjoint, after decreasing $\varepsilon$ the compact sets $C_j(\tau)$ are pairwise disjoint as well.

Suppose that $\Omega(\tau)$ has no edge. By Lemma~\ref{lem:extreme}, whose proof uses only \ref{S5}, $\partial\Omega(\tau)\subset\Sigma(\tau)$. By Lemma~\ref{lem:smooth}, $\partial\Omega(\tau)$ is a Jordan curve, so it is connected and lies in one of the disjoint compact sets $C_j(\tau)$; hence $\Omega(\tau)=\conv C_j(\tau)$.

If this happened for a sequence $\tau_m\downarrow\tau_0$, then, passing to a subsequence, $j$ would be fixed. Since $C_i(\tau)\to\gamma_i$ for every $i$, we have $\Sigma(\tau)\to\Gamma_0$, and by the continuity of the convex hull in the Hausdorff distance \cite[\S1.8]{Schneider14}, $\Omega(\tau_m)\to\conv\Gamma_0$ and $\conv C_j(\tau_m)\to\conv\gamma_j$. Hence $\conv\Gamma_0=\conv\gamma_j$, a contradiction.
\end{proof}

\begin{corollary}\label{cor:strictseveral}
If $\conv\gamma\ne\conv\Gamma_i$ for every component $\gamma$ of $\Gamma_i$, for $i=0$ or $i=1$, then $d<D(\rho_0,\rho_1)$.
\end{corollary}

\begin{proof}
By Lemma~\ref{lem:edgesnear}, applied at $\Gamma_0$ or, after reflection, at $\Gamma_1$, the hulls have edges on an interval of levels. Hence $\lambda>0$ there by \ref{S6}, and Remark~\ref{rem:strict} applies.
\end{proof}

For embedded surfaces the chords can be counted. Call a component of $\Sigma(\tau)$ \emph{exposed} if it meets $\partial\Omega(\tau)$, and let $e(\tau)$ be the number of exposed components.

\begin{lemma}\label{lem:exposed}
Let $M$ be embedded and let $\tau$ be a regular value with $e(\tau)\ge2$. Then $\Omega(\tau)$ has at least $e(\tau)$ chords whose endpoints lie on different components of $\Sigma(\tau)$. Consequently,
\[
2\pi\lambda(\tau)\ \ge\ \ell(\tau)\ \ge\ e(\tau)\,\delta(\tau),\qquad \delta(\tau)=\min_{i\ne j}\operatorname{dist}(C_i,C_j)>0,
\]
where $C_1,\dots,C_{N(\tau)}$ are the components of $\Sigma(\tau)$.
\end{lemma}

\begin{proof}
Since $M$ is embedded and $\tau$ is regular, the $C_j$ are disjoint Jordan curves, so $\delta(\tau)>0$. Let $J=\partial\Omega(\tau)$, a Jordan curve by Lemma~\ref{lem:smooth}, and let $K_j=J\cap C_j$. By Lemma~\ref{lem:extreme}, $J\setminus\bigcup_jK_j$ is the union of the gaps.

Call a gap \emph{mixed} if its endpoints lie on different components, and let $U_j$ be the union of $K_j$ with all gaps having both endpoints in $K_j$. Each $U_j$ is closed: a limit of such gaps shrinks to a point of the closed set $K_j$.

Let $m$ be the number of mixed gaps, and suppose $m<e(\tau)$.  If $m=0$, the connected set $J$ is the disjoint union of the finitely many nonempty closed sets $U_j$, so $e(\tau)=1$.
 If $m\ge1$, then $J$ minus the mixed gaps consists of at most $m$ closed arcs. Each of them lies in a single $U_j$ by connectedness, so at most $m$ components are exposed.

Both contradict $m<e(\tau)$, so $m\ge e(\tau)$. The closure of a mixed gap is a chord of length at least $\delta(\tau)$, and the chords are contained in edges with disjoint interiors. The bound on $\lambda$ is \ref{S6}.
\end{proof}

By Proposition~\ref{prop:Bineq}, the density of the absolutely continuous part of $B''$ therefore satisfies, at a.e.\ regular level $\tau$ with $e(\tau)\ge2$,
\[
B''\ \ge\ \frac{1+B'^2}{B-\frac{1}{2\pi}e(\tau)\,\delta(\tau)} .
\]
Thus many exposed components, or widely separated ones, force the hulls to be strongly convex in $\tau$. The hull does not see components which are not exposed; the count of all components of the sections is studied in \cite{T26notes} by different means.

\subsection*{Acknowledgments}
I am grateful to Vladimir Klyachin for useful discussions on minimal tubes in higher codimension.

\providecommand{\bysame}{\leavevmode\hbox to3em{\hrulefill}\thinspace}
\providecommand{\MR}{\relax\ifhmode\unskip\space\fi MR }
\providecommand{\MRhref}[2]{%
  \href{http://www.ams.org/mathscinet-getitem?mr=#1}{#2}
}
\providecommand{\href}[2]{#2}

\end{document}